\input ./style/arxiv-general.cfg
\documentclass[MSNbibl,number,citesort,seceqn,dvips]{arxbj}
\makeatletter
   \@ifpackageloaded{graphicx}{}{\usepackage{graphicx}}
\makeatother

\volume{22}
\issue{3}
\pubyear{2016}
\firstpage{1745}
\lastpage{1769}
\doi{10.3150/15-BEJ710}
\docsubty{FLA}

\makeatletter
\def\cal{\mathcal}
\newtheorem{Lem}{Lemma}[section]

\newtheorem{Theor}{Theorem}[section]

\newproclaim{defin}{Definition}

\newcommand{\R}{\mathbb{ R}}

\def\pmb{\bolds}

\newcommand{\Xb}{\mathbf{X}}

\newcommand{\Vb}{\mathbf{V}}
\newcommand{\eb}{{\mathbf{e}}}
\newcommand{\Ab}{{\mathbf{A}}}
\newcommand{\Bb}{{\mathbf{B}}}
\newcommand{\Ub}{{\mathbf{U}}}
\newcommand{\Wb}{{\mathbf{W}}}
\newcommand{\Yb}{{\mathbf{Y}}}

\newcommand{\thetab}{{\pmb\theta}}

\newcommand{\ny}{n\rightarrow\infty}

\makeatother

\begin{document}
\begin{frontmatter}

\title{On high-dimensional sign tests}
\runtitle{On high-dimensional sign tests}

\begin{aug}
\author[A]{\inits{D.}\fnms{Davy}~\snm{Paindaveine}\thanksref{A}\ead[label=e1]{dpaindav@ulb.ac.be}}
\and
\author[B]{\inits{T.}\fnms{Thomas}~\snm{Verdebout}\corref{}\thanksref{B}\ead[label=e2]{thomas.verdebout@univ-lille3.fr}}
\address[A]{ECARES and D\'epartement de Math\'ematiques, Universit\'e
Libre de Bruxelles (ULB), Avenue F.D. Roosevelt 50, CP114/04,
B-1050 Brussels, Belgium.
\printead{e1}}

\address[B]{INRIA, EQUIPPE, Universit\'e Lille 3 and D\'{e}partement
de Math\'{e}matique, Universit\'{e} libre de Bruxelles (ULB), Boulevard du
Triomphe, CP210, B-1050 Brussels, Belgium.\\
\printead{e2}}

\end{aug}

%
\received{\smonth{6} \syear{2014}}
%
\revised{\smonth{10} \syear{2014}}

%
\begin{abstract}
Sign tests are among the most successful procedures in multivariate
nonparametric statistics.
In this paper, we consider several testing problems in multivariate
analysis, directional
statistics and multivariate time series analysis, and we show that,
under appropriate
symmetry assumptions, the fixed-$p$ multivariate sign tests remain
valid in the
high-dimensional case. Remarkably, our asymptotic results are \emph{universal},
in the sense that, unlike in most previous works in high-dimensional
statistics, $p$ may
go to infinity in an arbitrary way as $n$ does. We conduct simulations
that (i) confirm our asymptotic results,
(ii) reveal that, even for relatively large $p$, chi-square critical
values are to be
favoured over the (asymptotically equivalent) Gaussian ones and (iii)
show that,
for testing i.i.d.-ness against serial dependence in the
high-dimensional case,
Portmanteau sign tests outperform their competitors in terms of
validity-robustness.
\end{abstract}

%
\begin{keyword}
\kwd{high-dimensional tests}
\kwd{Portmanteau tests}
\kwd{sign tests}
\kwd{universal asymptotics}
\end{keyword}
\end{frontmatter}

\section{Introduction}\label{sec1}

Sign procedures that discard the magnitude of the observations to
rather focus on their direction from a given location are among the
most popular nonparametric techniques. Multivariate sign tests, in
particular, have been extensively studied in the last decades.
Multivariate location was considered in Randles \cite{Ran89,Ran00}, M{\"
o}tt{\"o}nen and Oja \cite{mooj95} and Hallin and Paindaveine \cite{HP02}, whereas
problems on (normalized)
covariance or scatter matrices were considered in Tyler~\cite
{Tyl1987}, D{\"u}mbgen \cite
{Dum98}, Hallin and Paindaveine \cite{HP06} and Hallin, Paindaveine
and Verdebout \cite{HPV10}. Multivariate sign tests were
also developed, for example, for the problem of testing i.i.d.-ness
against serial dependence (see Paindaveine \cite{Pai2009}), or for
testing for
multivariate independence (see Taskinen, Kankainen and Oja \cite
{Tasetal2003} and Taskinen, Oja and Randles \cite
{Tasetal2005}). Most references above actually focus on so-called \emph
{spatial sign tests}, that is, on tests that are based on (a~possibly
standardized version of) the signs $\Ub_i=\Xb_i/\|\Xb_i\|$,
$i=1,\ldots,n$, obtained by projecting the $p$-variate observations
$\Xb_i$,
$i=1,\ldots,n$ on the unit sphere of $\R^p$. In the sequel, sign tests
will refer to spatial sign tests.

Multivariate sign tests enjoy many desirable properties.
First, they are robust, since they do not require stringent parametric
assumptions, nor any moment conditions. The above projection on the
unit sphere also guarantees that robustness holds with respect to
possible outliers (observations with large magnitudes). Second, for
fixed $p$, sign tests further enjoy uniformly high asymptotic
efficiency unless $p$ is small: for $p$-dimensional location and serial
problems, the \emph{lower bound} of the asymptotic relative
efficiencies of sign tests with respect to their classical Gaussian
competitors is $((p-1)/p)^2$, while for $p$-dimensional problems
involving normalized scatter matrices (testing for sphericity, PCA,
etc.), this lower bound is $(p/(p+2))^2$. This shows that, for
moderate-to-large $p$, using sign tests instead of classical Gaussian
tests may barely have any cost in terms of asymptotic efficiency. Even
better, in all problems above, there are no finite upper bounds, so
that the asymptotic efficiency gain of using sign tests may be
arbitrarily large. These remarkable asymptotic efficiency properties of
sign tests that may be puzzling at first sight are actually in line
with the fact that, as $p\to\infty$, the signs $\Ub_i$, $i=1,\ldots,n$
asymptotically contain all relevant information since data points
concentrate more and more on a common sphere (see, e.g., Hall, Marron
and Neeman \cite{Halletal2005}).

The asymptotic results in the previous paragraph, however, relate to
the fixed-$p$ \mbox{large-$n$} setup, hence are not directly
interpretable in an $(n,p)$-asymptotics framework. In this paper, we
therefore study the asymptotic null behaviour of several multivariate
sign tests in the high-dimensional setup. Actually we show that the
classical sign tests, based on their usual fixed-$p$ asymptotic
chi-square critical values, remain valid in the high-dimensional case.
In this sense, sign tests are robust to high-dimensionality (see
Section~\ref{sec2}). Beyond validating for the first time the use of
several sign tests in the high-dimensional setup, our results put the
emphasis on an interesting robustness property of sign tests, namely
the fact that they allow for \emph{universal} $(n,p)$-asymptotics, in
the sense that the respective null asymptotic distributions hold
whenever $\min(n,p)\to\infty$. In contrast, $(n,p)$-asymptotic results
in the literature usually restrict in a stringent way how $p$ may go to
infinity as a function of $n$ -- typically, it is imposed that $p/n\to
c$ for some $c$ belonging to a given convex set $C\subset[0,\infty)$
(most often, $C=[0,1)$ or $C=(1,+\infty)$). Some asymptotic
investigations cover all ($n,p$)-``regimes'', but different regimes
provide different asymptotic distributions, or lead to different test
statistics, which jeopardizes practical implementation. To the best of
our knowledge, thus, sign tests are the only tests that can be applied
without making (i)~strong restrictions on $(n,p)$ and (ii) severe
distributional or moment assumptions.


A huge amount of research has been dedicated in the last decade to
high-dimensional hypothesis testing. Location problems have been
investigated in, for example,
Srivastava and Fujikoshi~\cite{SriFuj2006},
and Srivastava and Kubokawa \cite{SriKub2013} (see also the references
therein),
while numerous papers have considered problems related to covariance or
scatter matrices; see, among many others, Ledoit and Wolf \cite
{LedWol2002}, Onatski, Moreira and Hallin \cite
{Onaetal2013} and Jiang and Yang \cite{JiaYan2013}. In this paper, we
study the
$(n,p)$-asymptotic null distribution of sign tests for various problems.

First, we tackle problems related with high-dimensional directional
data, which are more and more common, for example, in magnetic
resonance (Dryden \cite{Dry2005}) or gene-expression (Banerjee \textit
{et al.}
\cite
{banerjee2003generative}). In Section~\ref{sphersec}, we provide the
$(n,p)$-asymptotic null distribution of the Rayleigh~\cite{Ray1919} test
statistic that addresses the problem of testing uniformity on the unit
sphere. High-dimensional tests for this problem were recently proposed
in Cai and Jiang \cite{CaJia12} and Cai, Fan and Jiang \cite
{Caietal2013}. We show that the Rayleigh
test, unlike the latter competitors, is robust to high-dimensionality
in a universal way (in the sense explained above) and can therefore be
used for any $(n,p)$-regime. In the same section, we treat another
important problem that is standard in directional statistics, namely
the spherical location problem in the context of rotationally symmetric
distributions (see Mardia and Jupp \cite{MarJup2000} and Ley \textit
{et al.}
\cite{Leyetal2013}), and show
that the Paindaveine and Verdebout \cite{PaiVer2013} sign test is also
robust to high-dimensionality.

Then in Section~\ref{sersec}, we consider the problem of testing for
i.i.d.-ness against serial dependence in the general multivariate case,
which is arguably the most fundamental goodness-of-fit testing problem
in multivariate time series analysis. To the best of our knowledge, we
provide here the first high-dimensional result for this problem by
deriving, under appropriate symmetry assumptions, the universal
$(n,p)$-asymptotic null distribution of the Paindaveine \cite{Pai2009}
Portmanteau-type sign test. Finally, in Section~\ref{indsec}, we tackle
the problem of testing for multivariate independence and the problem of
testing for sphericity about a specified center.

In Section~\ref{simusec}, we conduct different simulations. In
Section~\ref{simusec1}, a Monte-Carlo study confirms that, when
properly standardized, sign test statistics are $(n,p)$-asymptotically
normal under the null. Yet, as we also show through simulations in
Section~\ref{simusec2}, the traditional chi-square critical values
better approximate the exact ones than their (actually, liberal)
Gaussian counterparts, even for relatively large $p$, hence should be
favoured for practical purposes. Finally, in Section~\ref{simusec3}, we
show that, unlike its classical competitors, the Portmanteau sign test
from Paindaveine \cite{Pai2009} has null rejection frequencies that
are robust to
high-dimensionality and heavy tails. We end the paper with the \hyperref
[app]{Appendix}
that contains proofs of technical results. For the sake of
completeness, the proofs for the independence and sphericity problems
are provided in the supplemental article Paindaveine and Verdebout
\cite{PaiVer2013c} that also
reports some additional simulation results.

\section{High-dimensional sign tests}
\label{sec2}

Consider some generic testing problem involving the null
hypothesis $\mathcal{H}_0$ and the alternative hypothesis $\mathcal
{H}_1$, to be addressed on the basis of $p$-variate observations $\Xb
_i$, $i=1,\ldots,n$. Let $Q_p^{(n)}$ be a test statistic for this problem,
that, for fixed $p$, is asymptotically $\chi^2_{d_p}$ under the null,
where $d_p\to\infty$ as $p\to\infty$ (all sign test statistics
considered in this paper meet this property). Then the corresponding
fixed-$p$ test, $\phi_p^{(n)}$ say, rejects the null at asymptotic
level $\alpha$ whenever
\[
Q_p^{(n)}>\chi^2_{d_p,1-\alpha},
\]
where $\chi^2_{d,1-\alpha}$ stands for the upper $\alpha$-quantile of
the chi-square distribution with $d$ degrees of freedom. Now, it is
well known that, if $Z_d$ is chi-square with $d$ degrees of freedom, then
$
(Z_d-d)/\sqrt{2d}
$
weakly converges to the standard normal distribution as $d\to\infty$.
Since we assumed that $d_p\to\infty$ as $p\to\infty$, it may then be
\emph{expected} that the test $\phi_{\mathcal{N},p}^{(n)}$ that
rejects the null whenever
\[
Q_{\mathcal{N},p}^{(n)}=\frac{Q_p^{(n)}-d_p}{\sqrt{2d_p}} > \Phi ^{-1}(1-
\alpha) 
\]
(throughout, $\Phi$ denotes the c.d.f. of the standard normal distribution)
has asymptotic level $\alpha$ under the null when both $n$ and $p$
converge to infinity. Clearly, the larger $p$, the closer the
tests $\phi_{p}^{(n)}$ and $\phi_{\mathcal{N},p}^{(n)}$, which may
then be considered
equivalent in the $(n,p)$-asymptotic framework.

Of course, establishing that the sequence of tests $\phi_{\mathcal
{N},p}^{(n)}$ -- hence
also the sequence of tests $\phi_{p}^{(n)}$ -- is valid in the
high-dimensional setup, that is, has asymptotic level $\alpha$ under
the null when both $n$ and $p$ go to infinity, requires a formal proof.
All the more so that the result does not always hold true. Some test
statistics indeed need to be appropriately corrected to obtain valid
$(n,p)$-asymptotic results, while some others do not (see, e.g., Ledoit
and Wolf \cite
{LedWol2002}). The test statistics that do not need be corrected can
be called \emph{high-dimensional (HD-)robust}. In the setup above, the
test statistic $Q_p^{(n)}$ is thus HD-robust if, under the null,
\[
\frac{Q_p^{(n)}-d_p}{\sqrt{2 d_p}}
\]
converges weakly to a standard Gaussian random variable as $n$ and $p$
go to infinity.

The main goal of this paper is to show that for the problems enumerated
in the \hyperref[sec1]{Introduction}, the traditional sign test
statistics are, under
appropriate symmetry conditions, \emph{universally} HD-robust, in the
sense that HD-robustness is achieved without imposing any constraint on
the way $p$ goes to infinity as $n$ does.

\subsection{Testing uniformity on the unit sphere}
\label{sphersec}

Let the random $p$-vectors $\Ub_1, \ldots, \Ub_n$ be mutually
independent and identically distributed, with a common distribution
that is supported on the unit sphere $\mathcal{S}^{p-1}=\{\mathbf
{x}\in\R^p\dvt\|
\mathbf{x}\|=\sqrt{\mathbf
{x}'\mathbf{x}}=1\}$ of $\R^p$. An important problem in
directional statistics consists in testing the null
hypothesis $\mathcal
{H}_0$ that this common distribution is the uniform on $\mathcal
{S}^{p-1}$. As already mentioned in the \hyperref[sec1]{Introduction},
this problem has
been recently considered in the high-dimensional case in Cai and Jiang
\cite{CaJia12}
and Cai, Fan and Jiang \cite{Caietal2013}. Arguably, the most
classical test for this
problem is the Rayleigh \cite{Ray1919} test that rejects the null for large
values of
%
%
\begin{equation}
\label{raylstat} R_p^{(n)} = \frac{p}{n} \sum
_{i,j=1}^n \Ub_i^{\prime}
\Ub_j;
\end{equation}
see, for example, Mardia and Jupp \cite{MarJup2000}, Section~6.3.
Under $\mathcal{H}_0$, $\Ub
_i$ has mean zero and covariance matrix $\frac{1}{p}\mathbf{I}_p$ (see
Lemma~\ref{LemmaU}), where $\mathbf{I}_\ell$ denotes the $\ell
\times
\ell$ identity matrix, so that the multivariate CLT readily implies
that, for any fixed $p$, the asymptotic null distribution
of $R_p^{(n)}$ is
$\chi^2_p$. Therefore, Rayleigh's test rejects the null, at asymptotic
level $\alpha$, whenever $R_p^{(n)}>\chi^2_{p,1-\alpha}$.

Obviously, applying Rayleigh's test of uniformity
is only possible when the sample size $n$ is large enough, compared
to $p$, to make the CLT approximation reasonable. We now derive the
asymptotic distribution of Rayleigh's test statistic $R_p^{(n)}$ in the
high-dimensional case when both $p=p_n$ and $n$ go to infinity. As
announced in the \hyperref[sec1]{Introduction}, our approach is \emph
{universal}, in the
sense that, unlike in most works on high-dimensional statistics, $p_n$
may go to infinity in a totally arbitrary way (the only restriction
being that both $p_n$ and $n$ go to infinity). Note in particular that
the asymptotic null distribution of the tests of uniformity proposed
in Cai, Fan and Jiang~\cite{Caietal2013}, which are based on
statistics of the form $\min_{i,j}\arccos(\Ub_i'\Ub_j)$, depends on
the $(n,p)$-regime considered.

Basically, we will show that the $(n,p)$-asymptotic distribution of
(the standardized version of) $R_p^{(n)}$ is universally standard
normal. To
do so, rewrite Rayleigh's statistic as
%
%
\begin{equation}
R_p^{(n)} = \frac{p}{n} \Biggl( n + \sum
_{1\leq i\neq j\leq n}^n \Ub_i^{\prime}
\Ub_j \Biggr) 
= p + \frac{2p}{n} \sum
_{1\leq i < j\leq n} \Ub_i^{\prime}
\Ub_j 
,
\end{equation}
and consider the standardized statistic
%
%
\begin{equation}
\label{RSt1} R_{\mathcal{N},p}^{(n)}= \frac{R_p^{(n)}-p}{\sqrt{2p}} =
\frac{\sqrt{2p}}{n} \sum_{1\leq i < j\leq n} \Ub_i^{\prime}
\Ub_j.
\end{equation}
As we recalled above, the fixed-$p$ asymptotic null distribution
of $R_p^{(n)}$ is $\chi^2_p$, hence has mean $p$ and variance $2p$, which
makes the standardization in (\ref{RSt1}) most natural. The main result
of this section is the following.

%
\begin{Theor}
\label{maintheorem}
Let $p_n$ be an arbitrary sequence of positive integers converging
to $+\infty$ as $\ny$.
Assume that $\Ub_{ni}$, $i=1,\ldots,n$, $n=1,2,\ldots,$ is a triangular
array such that for any $n$, the random $p_n$-vectors $\Ub_{ni}$,
$i=1,\ldots,n$ are \mbox{i.i.d.} uniform on $\mathcal{S}^{p_n-1}$. Then
%
%
\begin{equation}
\label{RSt} R_{\mathcal{N},p}^{(n)}= \frac{R_p^{(n)}-p_n}{\sqrt
{2p_n}} =
\frac{\sqrt{2p_n}}{n} \sum_{1\leq i < j\leq n} \Ub_{ni}^{\prime}
\Ub_{nj}
\end{equation}
converges in distribution to the standard normal as $n\to\infty$.
\end{Theor}

See the \hyperref[app]{Appendix} for the proof. As we now explain, we
are also able to
deal with another famous testing problem in directional statistics,
namely the spherical location problem. The relevant distributional
setup for this problem is the class of so-called \emph{rotationally
symmetric distributions} on $\mathcal{S}^{p-1}$; see Mardia and Jupp
\cite{MarJup2000}
or Ley \textit{et al.} \cite{Leyetal2013} for details. In the absolutely
continuous case,
this corresponds to the semiparametric class of densities (with respect
to the surface area measure on $\mathcal{S}^{p-1}$) of the form
%
%
\begin{equation}
\label{rotsym} \mathbf{u}\mapsto c_{p,f} f\bigl(\mathbf{u}^{\prime
}
\thetab\bigr), \qquad\mathbf{u} \in{\cal S}^{p-1},
\end{equation}
where $\thetab\in\mathcal{S}^{p-1}$ is a location parameter,
$c_{p,f}\ (>0)$ is a normalization constant and $f\dvtx[-1,1] \to\R^+$ is
some monotone increasing function. The spherical location problem
consists in testing that $\thetab$ is equal to some given
vector $\thetab_0\in\mathcal{S}^{p-1}$, on the basis of a random
sample $\Ub_1,\ldots,\Ub_n$ from~(\ref{rotsym}). Signed-rank tests for
this problem were recently proposed in Paindaveine and Verdebout \cite
{PaiVer2013} (while Ley \textit{et al.} \cite
{Leyetal2013} developed the corresponding estimators of $\thetab$). In
particular, the sign-based version of these tests rejects the null
hypothesis ${\cal H}_0\dvt\thetab= \thetab_0$ whenever
\[
R_p^{(n)}(\thetab_0) = \frac{p-1}{n} \sum
_{i,j=1}^n \mathbf{U}_{i}'(
\thetab_0) {\mathbf U}_{j}(\thetab_0) >
\chi^2_{p-1,1-\alpha} ,
\]
where
\[
\mathbf{U}_{i}(\thetab_0) = \frac{(\mathbf{I}_p- \thetab_0 \thetab
_0^{\prime}) \Ub_i}{\| (\mathbf{I}_p-
\thetab
_0 \thetab_0^{\prime}) \Ub_i \|} ,\qquad i=1,
\ldots, n
\]
is the multivariate sign of the projection of $\Ub_i$ onto the tangent
space to $\mathcal{S}^{p-1}$ at $\thetab_0$ (note that since the $\Ub
_i$'s have an absolutely continuous distribution on the sphere, the
corresponding $\mathbf{U}_{i}(\thetab_0)$'s are well-defined almost surely).

Irrespective of the underlying infinite-dimensional nuisance $f$,
the $\mathbf{U}_{i}(\thetab_0)$'s, under the null, are \mbox{i.i.d.} with
a common distribution that is uniform on the hypersphere ${\mathcal
S}^{p-2}_{\thetab_0}= \{ \mathbf{u}\in{\mathcal
S}^{p-1}\dvt\mathbf{u}^{\prime}\thetab
_0=0\}$. The following result then follows from Theorem~\ref{maintheorem}.

%
\begin{Theor}
\label{maintheoremrotsym}
Let $p_n$ be an arbitrary sequence of positive integers converging
to $+\infty$ as $\ny$.
Assume that $\Ub_{ni}$, $i=1,\ldots,n$, $n=1,2,\ldots,$ is a triangular
array such that for any $n$, the random $p_n$-vectors $\mathbf{U}_{ni}$,
$i=1,\ldots,n$, are \mbox{i.i.d.} with a common rotationally symmetric
density over $\mathcal{S}^{p_n-1}$. Then, letting ${\mathbf
U}_{ni}(\thetab_0)
= (\mathbf{I}_p- \thetab_0 \thetab_0^{\prime}) \Ub_{ni}/\|
(\mathbf{I}_p-
\thetab_0
\thetab_0^{\prime}) \Ub_{ni} \|$, we have that
\[
R_{\mathcal{N},p}^{(n)}(\thetab_0) = \frac{R_p^{(n)}(\thetab
_0)-(p_n-1)}{\sqrt{2(p_n-1)}} =
\frac{\sqrt{2(p_n-1)}}{n} \sum_{1\leq i < j\leq n} \mathbf{U}_{ni}'(
\thetab_0) \mathbf{U}_{nj}(\thetab_0)
\]
converges in distribution to the standard normal as $n\to\infty$.
\end{Theor}

Of course, the resulting (universal) test of spherical location rejects
the null ${\cal H}_0\dvt\thetab= \thetab_0$ at asymptotic level
$\alpha$
whenever $R_{\mathcal{N},p}^{(n)}(\thetab_0)$ exceeds the upper
$\alpha$-quantile $\Phi
^{-1}(1-\alpha)$ of the standard Gaussian distribution. At the same
level, the universal test of uniformity on $\mathcal{S}^{p-1}$ rejects
the null when $R_{\mathcal{N},p}^{(n)}$ exceeds $\Phi^{-1}(1-\alpha
)$. As discussed in the
beginning of Section~\ref{sec2}, one may alternatively perform the
original fixed-$p$ chi-square tests, since chi-square and Gaussian
critical values are asymptotically equivalent as $p\to\infty$. The
objective of Section~\ref{simusec2} is to compare both types of
critical values through simulations.

We end this section by stressing that Theorem~\ref{maintheorem} is also
relevant in the context of high-dimensional location testing, in a
framework where the $p$-variate observations $\Xb_i$, $i=1,\ldots,n$
have \emph{independent spherical directions about $\thetab\ (\in\R^p)$}.
Throughout the paper, this will mean that $(\Xb_1',\ldots,\Xb_n')'$ is
equal in distribution to $(\thetab'+R_1\Ub_1',\ldots,\thetab
'+R_n\Ub
_n')'$, where (i) the $\Ub_i$'s are \mbox{i.i.d.} uniform
over $\mathcal
{S}^{p-1}$ and (ii) the $R_i$'s are arbitrary random variables such
that $\mathrm{ P}[R_i=0]=0$ for all~$i$. In particular, the $R_i$'s may
fail to be mutually independent and/or may be dependent of the $\Ub
_i$'s. Also, parallel to the \emph{generalized spherical distributions}
from Frahm \cite{fra2004} and Frahm and Jaekel \cite
{frajae2007,frajae2010}, the $R_i$'s do
not need to be nonnegative. In the sequel, instead of ``\emph{the
$p$-variate random vectors $\Xb_i$, $i=1,\ldots,n$ have independent
spherical directions about the origin of $\R^p$}'', we will simply write
``\emph{the $p$-variate random vectors $\Xb_i$, $i=1,\ldots,n$ have
independent spherical directions}''.

Assume then that $\Xb_{ni}$, $i=1,\ldots,n$, $n=1,2,\ldots$ is a
triangular array such that for any $n$, the random $p_n$-vectors $\Xb
_{ni}$, $i=1,\ldots,n$ have independent spherical directions
about $\thetab\in\R^{p_n}$, and consider the (location) testing problem
associated with the null $\mathcal{H}_0\dvt\thetab=\mathbf{0}$ and the
alternative $\mathcal{H}_1\dvt\thetab\neq\mathbf{0}$. It readily follows
from Theorem~\ref{maintheorem} that the test rejecting the null whenever
\[
\frac{R_p^{(n)}(\Xb_1/\|\Xb_1\|,\ldots, \Xb_n/\|\Xb_n\|
)-p_n}{\sqrt{2p_n}} > \Phi^{-1}(1-\alpha)
\]
has $(n,p)$-asymptotic level $\alpha$, irrespective of the way $p_n$
goes to infinity with $n$. This settles the high-dimensional null
distribution of the so-called \emph{spatial sign location test}; see,
for example, Oja \cite{Oja2010}, Chapter~6.

One might consider that the assumption that the observations have
independent spherical directions is restrictive. We point out, however,
that it is less stringent than the assumptions that observations are
mutually independent with a common spherically symmetric distribution,
and that testing for sphericity is one of the most treated problems in
high-dimensional hypothesis testing (which would not be the case if
sphericity never holds). If the null of sphericity is not rejected,
then practitioners may resort to the location tests above (and to the
tests we propose in the next sections). Now, if observations fail to
have a spherically symmetric distribution (more generally, if they do
not have independent spherical directions), performing marginal
standardization may bring us closer to spherical symmetry or
independent spherical directions. A thorough investigation of the
impact of such a ``whitening'' step on the asymptotic null behaviour of
the tests we consider is however beyond the scope of this paper, and we
therefore leave this for future research.



\subsection{Testing for i.i.d.-ness against serial dependence}
\label{sersec}

In univariate time series analysis, the daily practice for location
models such as ARMA or ARIMA is deeply rooted in the so-called \emph
{Box and Jenkins methodology}; see, for example, Brockwell and Davis
\cite
{BroDav1991} for
details. An important role in this methodology is played by diagnostic
checking procedures, such as \emph{Portmanteau tests}, that address the
null that the residuals of the model at hand are not serially
correlated (one often speaks of the null hypothesis of \emph
{randomness}). These tests typically reject the null for large values
of $\sum_{h=1}^H (n-h)({r}(h))^2$, where ${r}(h)$ denotes the lag-$h$
sample autocorrelation in the residual series. If autocorrelations are
computed in the series of residual signs rather than in the series of
residuals themselves, one obtains the ``generalized runs tests''
of Dufour, Hallin and Mizera \cite
{Dufetal1998}, that are robust to heteroscedasticity (for $H=1$, these
tests reduce to the celebrated \emph{runs} test of randomness, which
justifies the terminology).

Diagnostic checking also belongs to daily practice of multivariate time
series analysis, where Portmanteau tests are based on sums of squared
norms of lag-$h$ autocorrelation matrices; see, for example, L\"{u}tkepohl
\cite{Lut2005}.
The corresponding sign tests, that can be seen as multivariate
(generalized) runs tests, were developed in Paindaveine \cite
{Pai2009}. To the best
of our knowledge, Portmanteau-type tests have not yet been studied in
the high-dimensional case.

Let $\Xb_1,\ldots,\Xb_n$ be random $p$-vectors and consider the problem
of testing the null hypothesis of randomness (white noise) versus the
alternative of serial dependence. This problem can be addressed by
considering the sign-based autocorrelation matrices
\[
\mathbf{r}(h) = \frac{p}{n-h} \sum_{t=h+1}^n
\Ub_t \Ub_{t-h}',\qquad h=1,\ldots,H,
\]
with $\Ub_t=\Xb_t/\|\Xb_t\|$, $t=1,\ldots,n$. More precisely, the
resulting (fixed-$p$) test is the Paindaveine~\cite{Pai2009} test,
that rejects the
null of randomness at asymptotic level $\alpha$ whenever
%
%
\begin{equation}
\label{Tchideux} T_p^{(n)} = \sum
_{h=1}^H (n-h) \bigl\| \mathbf{r}(h)\bigr \|
^2_\mathrm{ Fr} >\chi^2_{Hp^2,1-\alpha},
\end{equation}
where $\|\Ab\|_\mathrm{ Fr}=(\operatorname{ Trace}(\Ab\Ab'))^{1/2}$
is the Frobenius norm of $\Ab$. This test is a natural sign-based
multivariate extension of the univariate Portmanteau tests described
above. 

Following the same methodology as in Section~\ref{sphersec}, we will
study the universal $(n,p)$-asymptotic null behaviour of a standardized
version of $T_p^{(n)}$, under the assumption that the observations have
independent spherical directions. Since
\begin{eqnarray*}
(n-h) \bigl\| \mathbf{r}(h) \bigr\| ^2_\mathrm{ Fr}
&=& \frac{p^2}{n-h} \sum
_{s,t=h+1}^n \bigl(\Ub_{s-h}'
\Ub_{t-h} \bigr) \bigl(\Ub_s'
\Ub_t \bigr)
\\
&= & 
%
p^2
+ \frac{2p^2}{n-h} \sum_{h+1\leq s<t\leq n} \bigl(
\Ub_{s-h}' \Ub_{t-h} \bigr) \bigl(
\Ub_s' \Ub_t \bigr) ,
\end{eqnarray*}
we will consider the standardization of $T_p^{(n)}$ given by
\[
T_{\mathcal{N},p}^{(n)}= \frac{
T_p^{(n)}
-Hp^2}{\sqrt{2Hp^2}} = \frac{\sqrt{2p^2}}{\sqrt{H}} \sum
_{h=1}^H \frac{1}{n-h} \sum
_{h+1\leq s<t\leq n} \bigl(\Ub_{s-h}'
\Ub_{t-h} \bigr) \bigl(\Ub_s'
\Ub_t \bigr) .
\]
Again, this standardization is in line with the fixed-$p$/large-$n$
(chi-square) null asymptotic distribution of $T_p^{(n)}$ in (\ref{Tchideux}).

%
%
\begin{Theor}\label{maintheoremserial}
Let $p_n$ be an arbitrary sequence of positive integers converging
to $+\infty$ as $\ny$.
Assume that $\Xb_{nt}$, $t=1,\ldots,n$, $n=1,2,\ldots,$ is a triangular
array such that for any $n$, the random $p_n$-vectors $\Xb_{nt}$,
$t=1,\ldots,n$ have independent spherical directions.
%
Then, letting $\Ub_{nt}=\Xb_{nt}/\|\Xb_{nt}\|$ for any $n,t$, we
have that
%
%
\begin{equation}
\label{TSt} T_{\mathcal{N},p}^{(n)}= \frac{
T_p^{(n)}
-Hp_n^2}{\sqrt{2Hp_n^2}} =
\frac{\sqrt{2p_n^2}}{\sqrt{H}} \sum_{h=1}^H
\frac{1}{n-h} \sum_{h+1\leq s<t\leq n} \bigl(
\Ub_{n,s-h}' \Ub_{n,t-h} \bigr) \bigl(
\Ub_{ns}' \Ub_{nt} \bigr),
\end{equation}
converges in distribution to the standard normal as $n\to\infty$.
\end{Theor}

As a direct consequence of Theorem~\ref{maintheoremserial} above,
the Paindaveine \cite{Pai2009} test statistic is universally HD-robust
in the sense
described in the beginning of Section~\ref{sec2}. As we will show
through simulations in Section~\ref{simusec3}, this is not the case for
its classical competitors.

\subsection{Testing for multivariate independence, testing for sphericity}
\label{indsec}

Consider now the problem of testing that the $p$-variate marginal $\Xb$
and $q$-variate marginal $\Yb$ of the random vector $(\Xb', \Yb')'$ are
independent, on the basis of a random sample $(\Xb_1', \Yb
_1')',\ldots
,  (\Xb_n', \Yb_n')'$. A sign test for this problem was introduced in
Taskinen, Kankainen and Oja \cite{Tasetal2003}. The spherical version
of this test is based on sign
covariance matrices of the form
\[
\mathbf{C}_n = \frac{1}{n} \sum_{i=1}^n
\Ub_i \Vb_i',
\]
where $\Ub_i:= \Xb_i/ \| \Xb_i \|$ and $\Vb_i:= \Yb_i / \| \Yb_i
\|$
$i=1, \ldots, n$, are the multivariate signs associated with $\Xb_i$
and $\Yb_i$, respectively. More precisely, for fixed $p$ and $q$, the
resulting test rejects the null of multivariate independence at
asymptotic level $\alpha$ whenever
%
%
\begin{equation}
\label{Ojastat} I_{p,q}^{(n)} = n pq \|\mathbf{C}_n
\|^2_\mathrm{ Fr} = \frac{pq}{n} \sum
_{i,j=1}^n \bigl(\Ub_i'
\Ub_j \bigr) \bigl( \Vb_i'
\Vb_j \bigr)
\\
> \chi^2_{pq,1-\alpha}.
\end{equation}
Below, we obtain the universal asymptotic distribution of a
standardized version of $I_{p,q}^{(n)}$, when the $\Xb_i$'s and
the $\Yb_i$'s
have independent spherical directions.

Adopting the same approach as in the previous sections, decompose the
test statistic $I_{p,q}^{(n)}$ into
\[
I_{p,q}^{(n)} = pq + \frac{2 pq}{n} \sum
_{1 \leq i <j \leq n} \bigl(\Ub_i'
\Ub_j \bigr) \bigl(\Vb_i'
\Vb_j \bigr) ,
\]
and consider the standardized statistic
%
%
\begin{equation}
\label{ISt} I_{\mathcal{N},p,q}^{(n)}= \frac{I_{p,q}^{(n)}-pq}{\sqrt
{2pq}} =
\frac{\sqrt{2 pq}}{n} \sum_{1 \leq i <j \leq n} \bigl(
\Ub_i' \Ub_j \bigr) \bigl(
\Vb_i' \Vb_j \bigr) .
\end{equation}
We then have the following universal $(n,p)$-asymptotic normality result.

%
\begin{Theor}
\label{maintheoremindep}
Let $p_n$ and $q_n$ be arbitrary sequences of positive integers such
that $\operatorname{ max}(p_n,q_n)$ converges to $+\infty$ as $\ny$.
Assume that $(\Xb_{ni}', \Yb_{ni}')'$, $i=1,\ldots,n$, $n=1,2,\ldots,$
is a triangular array such that \textup{(i)} for any $n$, the $p_n$-variate
random vectors $\Xb_{ni}$, $i=1,\ldots,n$ and $q_n$-variate
marginals $\Yb_{ni}$, $i=1,\ldots,n$ have independent spherical
directions, and such that \textup{(ii)} $\Ub_{ni}=\Xb_{ni} / \| \Xb_{ni}\|$ and
$\Vb_{ni}=\Yb_{ni} / \| \Yb_{ni}\|$ are independent for all $i,n$.
Then
%
%
\begin{equation}
\label{Ist2} I_{\mathcal{N},p,q}^{(n)}= \frac{I_{p,q}^{(n)}-p_n
q_n}{\sqrt{2p_n q_n}} =
\frac{\sqrt{2p_n q_n}}{n} \sum_{1\leq i < j\leq n} \Ub_{ni}'
\Ub_{nj} \Vb_{ni}'\Vb_{nj}
\end{equation}
converges in distribution to the standard normal as $n\to\infty$.
\end{Theor}

For the sake of completeness, the result is proved in the supplemental
article Paindaveine and Verdebout \cite{PaiVer2013c}.
Note that asymptotic normality is obtained even when only one of the
dimensions $p_n$, $q_n$ goes to infinity. Testing for multivariate
independence is intimately related to testing block-diagonality of
covariance matrices (the sign test in (\ref{Ojastat}) rejects the null
when the off-diagonal blocks of an empirical sign-based covariance
matrix is too large, in Frobenius norm).

Finally, another classical problem in multivariate analysis that is
linked to covariance matrices (in this case, the null of interest if
that the covariance matrix is proportional to the identity matrix) is
the problem of testing for sphericity. Since the seminal paper Ledoit
and Wolf \cite
{LedWol2002}, this problem -- that consists in testing that the common
distribution of \mbox{i.i.d.} random $p$-vectors $\Xb_1,\ldots,\Xb_n$
is spherically symmetric
-- has been treated in many papers on high-dimensional inference; see,
for example, Chen, Zhang and Zhong \cite{Cheetal2010}, Jiang and Yang
\cite
{JiaYan2013} and Zou \textit{et al.} \cite{Zouetal2013}.
When testing for sphericity about a specified center (without loss of
generality, about the origin of $\R^p$), the natural fixed-$p$ sign
test of sphericity rejects the null at asymptotic level $\alpha$ whenever
%
%
\begin{equation}
\label{HPstat} S_p^{(n)} = \frac{p(p+2)}{2n} \sum
_{i,j=1}^n \biggl( \bigl(\Ub_i'
\Ub_j \bigr)^2-\frac
{1}{p} \biggr) >
\chi^2_{d_p,1-\alpha} ,
\end{equation}
with $\Ub_i=\Xb_i/\|\Xb_i\|$, $i=1,\ldots,n$, and $d_p=(p-1)(p+2)/2$;
see Hallin and Paindaveine \cite{HP06} and Sirki\"{a} \textit{et al.}
\cite
{Siretal2009}. Using the methodology proposed
above, it can be shown that, under the null, the universal
$(n,p)$-asymptotic distribution of a standardized version
of $S_p^{(n)}$ is
standard normal under extremely mild assumptions. More precisely, we
have the following result (see the supplemental article Paindaveine and
Verdebout \cite
{PaiVer2013c} for a proof).

%
\begin{Theor}
\label{maintheoremsphericity}
Let $p_n$ be an arbitrary sequence of positive integers converging
to $+\infty$ as $\ny$.
Assume that $\Xb_{ni}$, $i=1,\ldots,n$, $n=1,2,\ldots,$ is a triangular
array such that for any $n$, the random $p_n$-vectors $\Xb_{ni}$,
$i=1,\ldots,n$ have independent spherical directions.
Then, letting $\Ub_{ni}=\Xb_{ni}/\|\Xb_{ni}\|$ for any $n,i$, we
have that
%
%
\begin{equation}
\label{Sst} S_{\mathcal{N},p}^{(n)}= \frac{S_p^{(n)}-\ell
(p_n)}{\sqrt{2\ell(p_n)}} =
\frac{p_n\sqrt{p_n+2}}{n\sqrt{p_n-1}} \sum_{1\leq i<j\leq n}^n \biggl(
\bigl(\Ub_{ni}' \Ub_{nj} \bigr)^2-
\frac{1}{p_n} \biggr)
\end{equation}
converges in distribution to the standard normal as $n\to\infty$.
\end{Theor}

Performing the test in (\ref{HPstat}) on centered observations $\Xb
_i-\hat{\thetab}$, $i=1,\ldots,n$ of course provides a test of
sphericity about an unspecified center $\thetab$. Recently, Zou
\textit{et al.}
\cite
{Zouetal2013} showed that this test is not robust to high
dimensionality, and proposed a robustified test that allows $p_n$ to
increase to infinity at most as fast as $n^2$ (hence, this test does
not allow for universal $(n,p)$-asymptotics).



\section{Simulations}
\label{simusec}

In this section, we first conduct a Monte-Carlo study to check the
validity of our universal asymptotic results in Theorems \ref
{maintheorem} and \ref{maintheoremserial}.
Then we investigate whether or not, for practical purposes, the
resulting Gaussian critical values should be favoured over the
(asymptotically equivalent) chi-square ones. Finally, we show that the
classical competitors of the Paindaveine \cite{Pai2009} sign test are
based on
statistics that severely fail to be HD-robust.

\subsection{Checking universal asymptotics}
\label{simusec1}

For every $(n,p)\in\{4, 30, 200, 1 000\}^2$, we generated $M=10\,000$
independent random samples of the form
\[
\Xb_{i}^{(p)},\qquad {i}=1,\ldots, n,
\]
from the $p$-dimensional standard normal distribution. Then we
evaluated the standardized statistics $T_{\mathcal{N},p}^{(n)}$ (see
(\ref{TSt}))
on each of these $M$ samples, and we also computed the standardized
statistic $R_{\mathcal{N},p}^{(n)}$ (see (\ref{RSt}))
on the corresponding samples of unit vectors
\[
\Ub_i^{(p)}=\frac{\Xb_{i}^{(p)}}{\|\Xb_{i}^{(p)}\|},\qquad i=1,\ldots,n.
\]
%

Clearly, in each case, samples are generated from the respective null
model, so that, according to our asymptotic results, the resulting
empirical distributions of both standardized test statistics considered
should be close to the standard normal for \emph{virtually any}
combination of ``large'' $n$ and $p$ values (``virtually any'' here
translates the universal asymptotics). To assess this, Figures~2 and~3
in supplemental article Paindaveine and Verdebout \cite{PaiVer2013c}
provide, for each $(n,p)$, histograms of the $M=10\,000$ corresponding
values of $R_{\mathcal{N},p}^{(n)}$ and $T_{\mathcal{N},p}^{(n)}$,
respectively. Inspection of these figures reveals the following:

\begin{longlist}[(iii)]

\item[(i)] The empirical distributions of $R_{\mathcal{N},p}^{(n)}$
and $T_{\mathcal{N},p}^{(n)}$
are clearly compatible with Theorems \ref{maintheorem} and \ref
{maintheoremserial}. For both statistics, the Gaussian approximation is
valid for moderate to large values of $n$ and $p$, irrespective of the
ratio $p/n$. This confirms our universal asymptotic results.

\item[(ii)] For small\vspace*{-1pt} $n$ ($n=4$), $R_{\mathcal{N},p}^{(n)}$ seems to
be asymptotically
Gaussian when $p \to\infty$, which illustrates the fixed-$n$
asymptotic results from Chikuse \cite{Chi1991}. On the contrary,
$T_{\mathcal{N},p}^{(n)}
$ cannot
be well approximated by a Gaussian distribution for small $n$.

\item[(iii)] The empirical distributions of $R_{\mathcal{N},p}^{(n)}$
and $T_{\mathcal{N},p}^{(n)}$ are
approximately (standardized) chi-square distributions for small $p$ and
moderate-to-large $n$ (i.e., $p=4$ and $n\geq30$), which is consistent
with the classical fixed-$p$ $n$-asymptotic results.
\end{longlist}

For the sake of completeness, the results of a similar study for the
statistic $I_{\mathcal{N},p,q}^{(n)}$ (see (\ref{Ist2})) are reported
in the supplemental
article Paindaveine and Verdebout \cite{PaiVer2013c}.

%


\subsection{Comparing critical values}
\label{simusec2}

For the sake of clarity, we focus here on the problem of testing
uniformity on the unit sphere, that was considered in Section~\ref
{sphersec}. Our main result there (Theorem~\ref{maintheorem}) justifies
the use, in high dimensions, of two tests:
\begin{itemize}
\item The first test, $\phi_{\mathcal{N},p}^{(n)}$ say, is based on
\emph{Gaussian
critical values}, and rejects the null at asymptotic level $\alpha$
whenever $
R_{\mathcal{N},p}^{(n)}
> \Phi^{-1}(1-\alpha)$.
\item The second test, that we will denote as $\phi_p^{(n)}$, is the
standard fixed-$p$ sign test, based on \emph{chi-square critical
values}. This test, that rejects the null at asymptotic level $\alpha$ if
\[
R_p^{(n)}>\chi^2_{p,1-\alpha}\qquad
\mbox{equivalently, if } R_{\mathcal{N},p}^{(n)}>\frac{\chi
^2_{p,1-\alpha}-p}{\sqrt{2p}},
\]
is of course $(n,p)$-asymptotically equivalent to $\phi_{\mathcal
{N},p}^{(n)}$.
\end{itemize}

\begin{figure}

\includegraphics{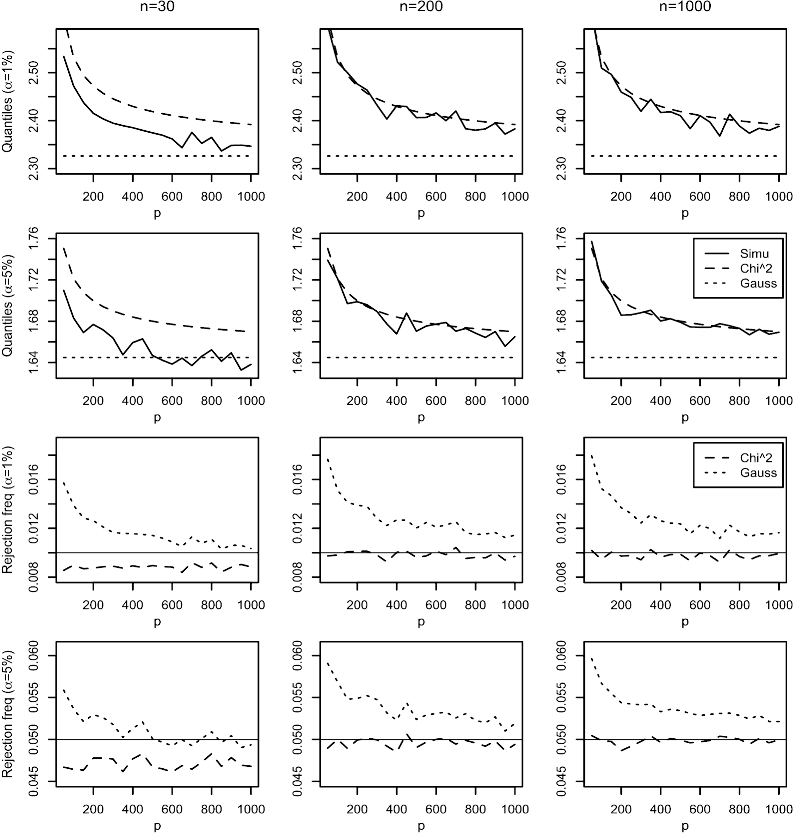}

\caption{For $\alpha=1\%$ and $\alpha=5\%$, the six upper panels
correspond to $n=30,200$ and $1 000$ and report the estimates $\hat
{q}_{n,p, 1-\alpha}$ of the exact upper $\alpha$-quantile (solid line)
of the statistic $R_{\mathcal{N},p}^{(n)}$, for $p=50,100,150,\ldots
,1 000$; the Gaussian
(dotted line) and chi-square (dashed line) approximations of these
exact quantiles are also provided. The six lower panels plot the
corresponding empirical rejection frequencies of
the Gaussian sign test $\phi_{\mathcal{N},p}^{(n)}$ (dotted line) and
chi-square sign
test $\phi_p^{(n)}$ (dashed line). Results are based on $M=100\,000$
independent replications; see Section~\protect\ref{simusec2} for details.}
\label{Rplot}
\end{figure}

To investigate what type of critical values should be favoured
depending on the $(n,p)$ configuration at hand, we performed the
following numerical exercise. In the exact same way as in the
simulations of Section~\ref{simusec1}, we generated, for each
$(n,p)\in
\{30, 200, 1 000\} \times\{50,100,150,\ldots,950,1 000\}$, $M=100\,000$ independent random samples
\[
\Ub_i^{(p)},\qquad i=1,\ldots,n
\]
from the uniform distribution over $\mathcal{S}^{p-1}$. For
each $(n,p)$ considered, the statistic $R_{\mathcal{N},p}^{(n)}$ was
evaluated on the
corresponding $M$ independent samples; for such a large $M$, the sample
upper $\alpha$-quantile, $\hat{q}_{n,p, 1-\alpha}$ say, of these $M$
independent replications of $R_{\mathcal{N},p}^{(n)}$ of course
provides a very accurate
estimate of the exact upper $\alpha$-quantile of this test statistic.
The appropriateness of Gaussian and chi-square critical values may thus
be evaluated by looking at how much these differ from $\hat{q}_{n,p, 1-
\alpha}$. This evaluation is made possible in the six upper panels of
Figure~\ref{Rplot}, which plot, for the various $n$ considered, the
sample quantiles $\hat{q}_{n,p, 0.95}$ and $\hat{q}_{n,p, 0.99}$ as a
function of $p$, along with the corresponding Gaussian and chi-square
critical values. Clearly, unless $n$ is small ($n=30$), chi-square
quantiles are to be favoured over Gaussian ones that tend to
underestimate the exact critical values.

The lower panels of Figure~\ref{Rplot} provide empirical rejection
frequencies of the Gaussian tests $\phi_{\mathcal{N},p}^{(n)}$ and of
the chi-square
tests $\phi_p^{(n)}$. In line with the results above, we conclude that,
when $n$ and $p$ are moderate to large (irrespective of the ratio
$p/n$), chi-square tests tend to achieve empirical type 1 risks that
are much closer to the nominal level than their Gaussian counterparts
that tend to be too liberal.

As a conclusion, while our universal asymptotic results broadly
validates the use of Gaussian -- hence, also of chi-square -- sign
tests, our recommendation is to avoid basing practical implementation
on Gaussian tests; instead, chi-square tests should be used in
practice. Finally, we point out that this conclusion is not only valid
for the problem of testing uniformity on the sphere, but does extend to
the other problems considered in this paper, as we checked by
performing similar numerical exercises -- with smaller $M$ values,
though, as the computational burden for the other tests is more severe
than for the Rayleigh test (this clearly makes $(n,p)$-asymptotic
critical values crucial for practical implementation).

%

\subsection{Comparing Portmanteau tests}
\label{simusec3}

In this last simulation exercise, we compare the high-dimensional
behaviour of the Portmanteau sign test described in Section~\ref
{sersec} with those of some classical (fixed-$p$) competitors, namely
the Ljung and Box \cite{ljung78} test (LB) and the Li and McLeod \cite
{li1981} test (LM). To the
best of our knowledge, no tests are currently available in the
high-dimensional case for this problem.

%
\begin{figure}

\includegraphics{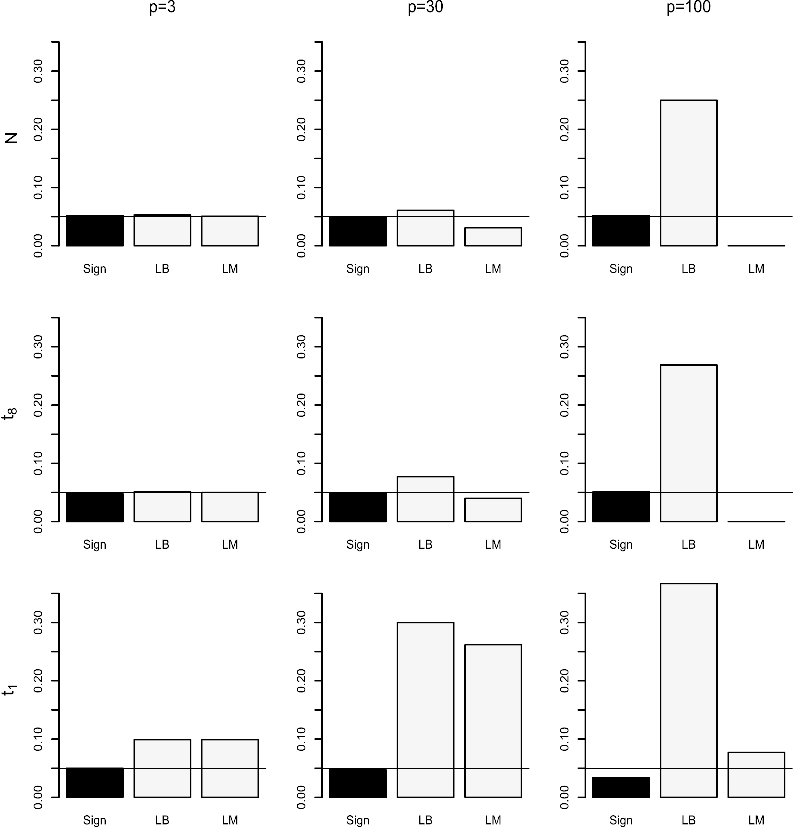}

\caption{\small Empirical rejection frequencies of the chi-square
Paindaveine \cite
{Pai2009} sign test (in black) and of the LM and LB tests (in light
grey). All tests were performed at the nominal level $5\%$. Results are
based on $M=1 000$ replications and the sample size is throughout $n=150$.}
\label{plotrunscomp}
\end{figure}

In order to do so, we generated, for every $p \in\{3, 30, 100\}$, $M=1
000$ independent random samples of the form
\[
\Xb_{jt}^{(p)},\qquad j=1,2, 3, {t}=1,\ldots, n=150.
\]
%
The $\Xb_{1t}^{(p)}$'s are standard Gaussian, whereas the $\Xb
_{2t}^{(p)}$'s (resp., the $\Xb_{3t}^{(p)}$'s) are (standard) student
with $8$ degrees of freedom (resp., with one degree of freedom).
For all of these (null) samples, we performed the three Portmanteau
tests mentioned above at nominal level $\alpha=5\%$. The resulting
rejection frequencies are represented in Figure~\ref{plotrunscomp}.

Inspection of this figure reveals that the Paindaveine \cite{Pai2009}
sign test is
the only test that is robust to high-dimensionality (and to heavy
tails). In the high-dimensional case, the LB test is extremely liberal,
while the LM test is liberal under the Cauchy but basically never
rejects the null under the Gaussian and the $t_8$. It should be noted
that the reason why we restricted here to cases for which $n>p$ is that
the LB and LM tests collapse when $n<p$ (in the sense that the
corresponding test statistics then cannot even be computed). In
contrast, the simulations of Section~\ref{simusec1} showed that the
sign test perfectly can deal with such cases.


%
\begin{appendix}\label{app}

\section{Some preliminary lemmas}
\label{appsec}


%
\begin{Lem}
\label{Lemmamoments}
Let $\Ub_1, \ldots, \Ub_n$ be \mbox{i.i.d.} uniform on ${\cal
S}^{p-1}$, and write $\rho_{ij}=\Ub_i'\Ub_j$. Then, for any $i,j$, \textup{(i)}
$\rho_{ij}^2\sim\operatorname{ Beta}(1/2,(p-1)/2)$, that is, $\rho_{ij}^2$
follows the Beta distribution with parameters $1/2$ and $(p-1)/2$; \textup{(ii)}
for any odd positive integer $m$, $\mathrm{ E}[\rho_{ij}^m]=0$; \textup{(iii)} for
any even positive integer $m$,
\[
\mathrm{ E} \bigl[\rho_{ij}^m \bigr] =\prod
_{r=0}^{m/2} \biggl( \frac{1+2r}{p+2r} \biggr),
\]
so that
\[
\mathrm{ E} \bigl[\rho_{ij}^2 \bigr]=\frac{1}{p},\qquad
\mathrm{ E} \bigl[\rho_{ij}^4 \bigr]=\frac{3}{p(p+2)},\qquad
\mathrm{ E} \bigl[\rho_{ij}^6 \bigr]=\frac{15}{p(p+2)(p+4)}
\]
and
\[
\mathrm{ E} \bigl[\rho_{ij}^8 \bigr]=\frac{105}{p(p+2)(p+4)(p+6)};
\]
\textup{(iv)} the $\rho_{ij}$'s, $i<j$, are pairwise independent (hence
uncorrelated); \textup{(v)} fix $h\in\{1,2,\ldots,n-2\}$. Then, for
any $i,j,s,t\in\{h+1,\ldots,n\}$ with $i<j$ and $s<t$, $\rho
_{i-h,j-h}\rho_{i,j}$ and $\rho_{s-h,t-h}\rho_{s,t}$ are uncorrelated,
unless $(i,j)=(s,t)$, in which case $\operatorname{ Cov}[\rho_{i-h,j-h}\rho
_{i,j},\break \rho_{s-h,t-h}\rho_{s,t}]$
$=1/p^2$.
\end{Lem}

\begin{pf}
(i) Rotational invariance of
the uniform distribution on $\mathcal{S}^{p-1}$ readily implies
that $\rho_{ij}$ is equal in distribution to $\eb_{p,1}'\Ub_j$,
where $\mathbf{e}_{p,\ell}$ denotes the $\ell$th vector in the canonical
basis of $\R^p$. The result then follows from the fact that $(\eb
_1'\Ub
_j)^2\sim\operatorname{ Beta}(1/2,(p-1)/2)$; see, for example, Muirhead
\cite{Mui1982},
Theorem~1.5.7(ii).
(ii) This is a trivial corollary of the fact that $\Ub_i$ and $-\Ub_i$,
hence also $\rho_{ij}$ and $-\rho_{ij}$, are equal in distribution.
(iii) This directly follows from (i) and the fact that, if $Y\sim
\operatorname{
Beta}(\alpha,\beta)$, then
\[
\mathrm{ E} \bigl[Y^s \bigr] = \prod_{r=0}^{s-1}
\biggl( \frac{\alpha+r}{\alpha+\beta+r} \biggr),
\]
for any positive integer $s$; see, for example, Johnson, Kotz and
Balakrishnan \cite{Johetal1995},
equation (25.14). (iv) If $i,j,r,s\in\{1,2,\ldots,n\}$, are pairwise
different, $\rho_{ij}$ and $\rho_{rs}$ are trivially independent. Let
then $i,j,s$ be three different integers in $\{1,2,\ldots,n\}$. Then
the rotational invariance of the uniform distribution on $\mathcal
{S}^{p-1}$ entails that the joint distribution of $(\rho_{ij},\rho
_{ir})$ coincides with that of $(\eb_{p,1}'\Ub_j,\eb_{p,1}'\Ub_r)$,
which has independent marginals.\vspace*{1pt} The result follows. (v) From
parts (iv) and (ii) of the lemma, we obtain
%
%
\begin{equation}
\label{expecov} \operatorname{ Cov}[\rho_{i-h,j-h}\rho_{i,j},
\rho_{s-h,t-h}\rho_{s,t}]=\mathrm{ E}[\rho_{i-h,j-h}
\rho_{s-h,t-h}\rho_{i,j}\rho_{s,t}].
\end{equation}
If $j\neq t$, this expectation is equal to zero, since
\[
(\Ub_1,\ldots, \Ub_{\max(j,t)-1},\pm\Ub_{\max(j,t)},
\Ub_{\max(j,t)+1},\ldots, \Ub_{n})
\]
are equal in distribution. Similarly, using the fact that
\[
(\Ub_1,\ldots, \Ub_{\min(i-h,s-h)-1},\pm\Ub_{\min(i-h,s-h)}, \Ub
_{\min
(i-h,s-h)+1},\ldots, \Ub_{n})
\]
are equal in distribution, we obtain that the expectation in (\ref
{expecov}) is equal to zero. Thus, to obtain a nonzero covariance, we
need to have $(i,j)=(s,t)$, which leads to
\[
\operatorname{ Cov}[\rho_{i-h,j-h}\rho_{i,j},\rho_{s-h,t-h}
\rho_{s,t}]=\operatorname{ Var}[\rho_{s-h,t-h}\rho_{s,t}]=
\mathrm{ E} \bigl[\rho_{s-h,t-h}^2 \bigr]\mathrm{ E} \bigl[\rho
_{s,t}^2 \bigr]=1/p^2.
\]
\upqed\end{pf}

Let $\mathbf{J}_p= \sum_{i,j =1}^{p}(\mathbf{e}_{p,i}\mathbf
{e}_{p,j}^{\prime}
)\otimes
(\mathbf{e}_{p,i}\mathbf{e}_{p,j}^{\prime}) = (\operatorname{ vec}
\mathbf{I}_p)
(\operatorname{ vec}
\mathbf{I}_p)^{\prime}$ and
consider the \emph{commutation matrix}
$\mathbf{K}_{p}= \sum_{i,j =1}^{p}(\mathbf{e}_{p,i}\mathbf
{e}_{p,j}^{\prime})\otimes
(\mathbf{e}_{p,j}\mathbf{e}_{p,i}^{\prime})$.
We then have the following result.

%
\begin{Lem}
\label{LemmaU}
Let $\Ub,\Vb,\Wb$ be mutually independent and uniformly distributed
on ${\cal S}^{p-1}$, ${\cal S}^{q-1}$, and ${\cal S}^{m-1}$,
respectively. Then \textup{(i)} $\mathrm{ E}[\Ub]=\mathbf{0}$,
\textup{(ii)} $\operatorname
{ Var}[\Ub
]=\frac{1}{p}\mathbf{I}_{p}$, \textup{(iii)}
$
\mathrm{ E} [\operatorname{ vec}(\Ub\Ub') (\operatorname{ vec}(\Ub\Ub') )'
]
=
\frac{1}{p(p+2)} ( \mathbf{I}_{p^2} + \mathbf{K}_p + \mathbf
{J}_p )$,
\textup{(iv)}
$
\mathrm{ E} [\operatorname{ vec}(\Ub\Vb') (\operatorname{ vec}(\Ub\Vb') )'
]
=
\frac{1}{pq} \mathbf{I}_{pq}$,
and
\textup{(v)}
$
\mathrm{ E} [ \operatorname{vec}(\Ub\Vb') \times (\operatorname{ vec}(\Ub\Wb') )'
]
=
\mathbf{0}_{pq\times pm}$, where $\mathbf{0}_{k\times\ell}$ denotes
the $k \times\ell$ zero matrix.
\end{Lem}

\begin{pf}
(i)--(ii) These identities follow directly from the orthogonal
invariance of $\Ub$ and the fact that $\|\Ub\|=1$ almost surely.
(iii) See Tyler \cite{Tyl1987}, page 244.
(iv)
The independence of $\Ub=(U_1,\ldots,U_p)'$ and $\Vb=(V_1,\ldots,V_q)'$
readily gives
\[
\mathrm{ E} \bigl[\operatorname{ vec} \bigl(\Ub\Vb' \bigr) \bigl(
\operatorname{ vec} \bigl(\Ub\Vb' \bigr) \bigr)' \bigr] = {
\sum_{i,j=1}^p \sum
_{r,s=1}^q} \mathrm{ E}[U_iU_j]
\mathrm{ E}[V_rV_s] \operatorname{ vec} \bigl(
\eb_{p,i}\eb_{q,r}' \bigr) \bigl(\operatorname{ vec}
\bigl(\eb_{q,j}\eb_{q,s}' \bigr)
\bigr)' ,
\]
which yields the result since $\mathrm{ E}[U_iU_j]=\frac{1}{p}\delta_{ij}$
and $\mathrm{ E}[V_rV_s]=\frac{1}{q}\delta_{rs}$ (see (ii)).
(v) This directly follows from the fact that $\pm\Wb$ are equal in
distribution.
\end{pf}

%
\begin{Lem}
\label{Lemmaseries}
For any $a>0$, $\sum_{\ell=1}^{n} \ell^a=\mathrm{O}(n^{a+1})$ as
$n\to\infty$.
\end{Lem}

\begin{pf}
The result follows from the fact that the sequence $(c_n)$, defined by
$
c_n=\frac{1}{n} \sum_{\ell=1}^{n} (\frac{\ell}{n} )^a
$,
is a sequence of Riemann sums for $\int_0^1 x^a\, \mathrm{d}x$, hence
converges
in $\R$.
\end{pf}

\section{Proofs of Theorems \texorpdfstring{\protect\ref{maintheorem}}{2.1} and \texorpdfstring{\protect\ref{maintheoremserial}}{2.3}}

{\sc Proof of Theorem~\ref{maintheorem}}. To prove this result (and the
corresponding results in the next two sections), we adopt an approach
that exploits the martingale difference structure of some process. In
that framework, the key result we will use is Theorem~35.12
from Billingsley \cite
{Bil1995}. Since this result plays a crucial role in the paper, we
state it here, in a form that is suitable for our purposes.

%
\begin{Theor}[(Billingsley \cite{Bil1995}, Theorem~35.12)]
\label{Billingsley}
Let $D_{n\ell}$, $\ell=1,\ldots,n$, $n=1,2,\ldots,$ be a triangular
array of random variables such that, for any $n$, $D_{n1},D_{n2},\ldots
,D_{nn}$ is a martingale difference sequence with respect to some
filtration $\mathcal{F}_{n1},\mathcal{F}_{n2},\ldots,\mathcal{F}_{nn}$.
Assume that, for any $n,\ell$, $D_{n\ell}$ 
has a finite variance. Letting $\sigma^2_{n\ell}=\mathrm{ E}
[D_{n\ell
}^2 | \mathcal{F}_{n,\ell-1} ]$ (with $\mathcal{F}_{n0}$ being
the trivial $\sigma$-algebra $\{\varnothing,\Omega\}$ for all $n$),
further assume that, as $n\to\infty$,
%
%
\begin{equation}
\label{Condition1} \sum_{\ell=1}^n
\sigma^2_{n\ell} \stackrel{\mathrm{ P}} {\to} 1
\end{equation}
(where $\stackrel{\mathrm{ P}}{\to}$ denotes convergence in
probability), and
%
%
\begin{equation}
\label{Condition2} \sum_{\ell=1}^n \mathrm{ E}
\bigl[D_{n\ell}^2 \mathbb{I}\bigl[|D_{n\ell
}|>\varepsilon\bigr]
\bigr]\to0.
\end{equation}
Then $\sum_{\ell=1}^n D_{n\ell}$ is asymptotically standard normal.
\end{Theor}

In order to apply this result, we define ${\cal F}_{n\ell}$ as the
$\sigma$-algebra generated by $\Ub_{n1}, \ldots, \Ub_{n\ell}$, and,
writing $\mathrm{ E}_{n\ell}$ for the conditional expectation with respect
to ${\cal F}_{n\ell}$, we let
\[
D^R_{n\ell} = \mathrm{ E}_{n\ell}
\bigl[R_{\mathcal{N},p}^{(n)}\bigr]-\mathrm{ E}_{n,\ell-1}
\bigl[R_{\mathcal{N},p}^{(n)}\bigr] = \frac{\sqrt{2p_n}}{n} \sum
_{i=1}^{\ell-1} \Ub_{ni}^{\prime}
\Ub_{n\ell} ,
\]
for any $\ell=1,2,\ldots$ (throughout, sums over empty set of indices
are defined as being equal to zero). Clearly, we have that
\[
R_{\mathcal{N},p}^{(n)}=\frac{\sqrt{2p_n}}{n} \sum
_{1\leq i<j\leq
n} { \Ub}_{ni}^{\prime}\Ub_{nj}
= \sum_{\ell=1}^n D^R_{n\ell}.
\]
Since $|D^R_{n\ell}|\leq\sqrt{2p_n}(\ell-1)/n$ almost surely,
every $D^R_{n\ell}$ has a finite-variance. Therefore, to establish
Theorem~\ref{maintheorem}, it is sufficient to prove the following two
lemmas, that show that (\ref{Condition1})--(\ref{Condition2}) are
fulfilled in the present context.

%
\begin{Lem}
\label{firstlemspher}
Letting $\sigma^2_{n\ell}=\mathrm{ E}_{n,\ell-1} [(D^R_{n\ell
})^2
]$, $\sum_{\ell=1}^n \sigma^2_{n\ell}$ converges to one in quadratic
mean as $n\to\infty$.
\end{Lem}

%
\begin{Lem}
\label{THElemmaspher}
\mbox{ For any $\varepsilon>0$, $\sum_{\ell=1}^n \mathrm{ E}
[(D^R_{n\ell})^2 {\mathbb I}[| D^R_{n\ell}| > \varepsilon] ] \to
0$ as $n\to\infty$.}
\end{Lem}

The proofs of these lemmas make intensive use of the properties of the
inner products $\rho_{n,ij}={\Ub}_{ni}^{\prime}\Ub_{nj}$; see
Lemma~\ref
{Lemmamoments} in Appendix~\ref{appsec}. For notational simplicity, we
will systematically drop the dependence on $n$ in $\rho_{n,ij}$, $\Ub
_{ni}$, $\mathrm{ E}_{n\ell}$, $\mathrm{ E}_n$ and $\operatorname{Var}_n$; here, $\mathrm{
E}_n$ and $\operatorname{ Var}_n$ stand for the unconditional expectation and
unconditional variance computed with respect to the joint distribution
of the $\Ub_{ni}$'s, $i=1,\ldots,n$.

\begin{pf*}{Proof of Lemma~\ref{firstlemspher}}
Using Lemma~\ref{LemmaU}(i)--(ii), we obtain
\[
\sigma^2_{n\ell} = \mathrm{
E}_{\ell-1} \bigl[ \bigl(D^R_{n\ell}
\bigr)^2 \bigr] = \frac{2p_n}{n^2} \sum
_{i,j=1}^{\ell-1} \Ub_i^{\prime}\mathrm{
E}\bigl[ \Ub_\ell\Ub_\ell^{\prime}\bigr]
\Ub_j = \frac{2}{n^2} \sum_{i,j=1}^{\ell-1}
\rho_{ij} , 
\]
which yields (Lemma~\ref{Lemmamoments}(ii))
%
%
\begin{equation}
\label{tobeusedspher} \mathrm{ E} \bigl[\sigma^2_{n\ell} \bigr] =
\frac{2(\ell-1)}{n^2} \qquad\bigl( = \mathrm{ E} \bigl[ \bigl(D^R_{n\ell}
\bigr)^2 \bigr] \bigr) .
\end{equation}
Therefore, as $n\to\infty$,
\[
\mathrm{ E} \Biggl[\sum_{\ell=1}^n
\sigma^2_{n\ell} \Biggr] = \frac{2}{n^2} \sum
_{\ell=1}^n (\ell-1) = \frac{n-1}{n} \to1.
\]
Using Lemma~\ref{Lemmamoments}(iv), then the fact that $\operatorname{
Var}[\rho
_{ij}]=1/p_n$ (which follows from Lemma~\ref{Lemmamoments}(ii)--(iii)),
we obtain that
\begin{eqnarray*}
\operatorname{ Var} \Biggl[\sum_{\ell=1}^n
\sigma^2_{n\ell} \Biggr] &= &\frac{16}{n^4} \operatorname{ Var}
\Biggl[ \sum_{\ell=3}^n \sum
_{1\leq i<j\leq\ell-1} \rho_{ij} \Biggr] = \frac{16}{n^4}
\operatorname{ Var} \biggl[ \sum_{1\leq i<j\leq n} (n-j)
\rho_{ij} \biggr]
\\
&  = &\frac{16}{n^4 p_n} \sum_{1\leq i<j\leq n}
(n-j)^2 = \frac{16}{n^4 p_n} \sum_{j=2}^n
(j-1) (n-j)^2
\\
&  =& \frac{16}{n^4 p_n} \sum_{j=1}^{n-1} j
(n-j-1)^2 \leq\frac{16}{n^2 p_n} \sum_{j=1}^{n-1}
j = \frac{8(n-1)}{n p_n},
\end{eqnarray*}
which is $\mathrm{o}(1)$ as $n\to\infty$. The result follows.
\end{pf*}

\begin{pf*}{Proof of Lemma~\ref{THElemmaspher}}
Applying first the Cauchy--Schwarz inequality, then the Chebyshev
inequality (note that $D^R_{n\ell}$ has zero mean), yields
\begin{eqnarray*}
\sum_{\ell=1}^n \mathrm{ E} \bigl[
\bigl(D^R_{n\ell} \bigr)^2 {\mathbb I} \bigl[\bigl|
D^R_{n\ell
}\bigr| > \varepsilon \bigr] \bigr] &\leq& \sum
_{\ell=1}^n \sqrt{\mathrm{ E} \bigl[
\bigl(D^R_{n\ell} \bigr)^4 \bigr]} \sqrt
{ \mathrm{ P} \bigl[\bigl|D^R_{n\ell}\bigr| > \varepsilon
\bigr]}
\\
&\leq& \frac{1}{\varepsilon} \sum_{\ell=1}^n
\sqrt{\mathrm{ E} \bigl[ \bigl(D^R_{n\ell}
\bigr)^4 \bigr]} \sqrt{\operatorname{ Var} \bigl[D^R_{n\ell}
\bigr]} .
\end{eqnarray*}
From (\ref{tobeusedspher}), we readily obtain that
$
\operatorname{ Var}[D^R_{n\ell}]
\leq
2(\ell-1)/n^2
$, which provides
%
%
\begin{equation}
\sum_{\ell=1}^n \mathrm{ E} \bigl[
\bigl(D^R_{n\ell} \bigr)^2 {\mathbb I} \bigl[\bigl|
D^R_{n\ell
}\bigr| > \varepsilon \bigr] \bigr] \leq
\frac{\sqrt{2}}{\varepsilon n} \sum_{\ell=1}^n \sqrt
{(\ell-1) \mathrm{ E} \bigl[ \bigl(D^R_{n\ell}
\bigr)^4 \bigr]} . \label{mainspher}
\end{equation}
Now, Lemma~\ref{Lemmamoments}(iv) yields
\begin{eqnarray*}
\mathrm{ E} \bigl[ \bigl(D^R_{n\ell} \bigr)^4
\bigr] &= &\frac{4p_n^2}{n^4} \mathrm{ E} \Biggl[ \Biggl( \sum
_{i=1}^{\ell-1} \rho_{i\ell} \Biggr)^4
\Biggr] = \frac{4p_n^2}{n^4} \sum_{i,j,r,s=1}^{\ell-1}
\mathrm{ E} [ \rho_{i\ell}\rho_{j\ell}\rho_{r\ell}
\rho_{s\ell} ]
\nonumber
\\[-8pt]
\\[-8pt]
\nonumber
&  =& \frac{4p_n^2}{n^4} \bigl\{ (\ell-1) \mathrm{ E} \bigl[
\rho_{i\ell}^4 \bigr] + 3(\ell-1) (\ell-2) \mathrm{ E} \bigl[
\rho_{1\ell}^2 \bigr] \mathrm{ E} \bigl[\rho_{2\ell}^2
\bigr] \bigr\},
\end{eqnarray*}
so that, using Lemma~\ref{Lemmamoments}(iii), we obtain
%
%
\begin{equation}
\mathrm{ E} \bigl[ \bigl(D^R_{n\ell} \bigr)^4
\bigr] = \frac{4p_n^2}{n^4} \biggl\{ \frac{3(\ell-1)}{p_n(p_n+2)} + \frac{3(\ell-1)(\ell-2)}{p_n^2}
\biggr\} \leq\frac{24(\ell-1)^2}{n^4} \label{okspher} .
\end{equation}
Plugging (\ref{okspher}) into (\ref{mainspher}), we conclude that
\begin{eqnarray*}
\sum_{\ell=1}^n \mathrm{ E} \bigl[
\bigl(D^R_{n\ell} \bigr)^2 {\mathbb I} \bigl[\bigl|
D^R_{n\ell
}\bigr| > \varepsilon \bigr] \bigr]
&\leq&
\frac{\sqrt{2}}{\varepsilon n} \sum_{\ell=1}^n \sqrt{
\frac{24(\ell-1)^3}{n^4}} \\
&\leq&\frac{\sqrt{48}}{\varepsilon n^3} \sum_{\ell=1}^n
(\ell-1)^{3/2} ,
\end{eqnarray*}
which is $\mathrm{O}(n^{-1/2})$ (see Lemma~\ref{Lemmaseries}). The
result follows.
\end{pf*}

\begin{pf*}{Proof of Theorem~\ref{maintheoremserial}}
We define ${\cal
F}_{n\ell}$ as the $\sigma$-algebra generated by $\Xb_{n1}, \ldots,
\Xb
_{n\ell}$, and we let
\[
D^T_{n\ell} = \mathrm{ E}_{n\ell}
\bigl[T_{\mathcal{N},p}^{(n)}\bigr]-\mathrm{ E}_{n,\ell-1}
\bigl[T_{\mathcal{N},p}^{(n)}\bigr] = \frac{\sqrt
{2p_n^2}}{\sqrt{H}} \sum
_{h=1}^H \frac{1}{n-h} \sum
_{s=h+1}^{\ell-1} \rho_{n,s-h,\ell-h} \rho_{n,s\ell}
\]
(recall that sums over empty sets of indices are defined as zero),
where we wrote $\rho_{n,s\ell}=\Ub_{ns}' \Ub_{nt}$ and where
$\mathrm{
E}_{n\ell}$ still denotes conditional expectation with respect
to ${\cal F}_{n\ell}$. This provides $
T_{\mathcal{N},p}^{(n)}
=
\sum_{\ell=1}^n D^T_{n\ell}$,
where $D^T_{n\ell}$ is almost surely bounded, hence has a
finite-variance. As in the previous section, asymptotic normality is
then proved by applying Theorem~\ref{Billingsley}, which is based on
both following lemmas.
\end{pf*}

%
\begin{Lem}
\label{firstlemserial}
Letting $\sigma^2_{n\ell}=\mathrm{ E}_{\ell-1} [(D^T_{n\ell
})^2 ]$,
$\sum_{\ell=1}^n \sigma^2_{n\ell}$ converges to one in quadratic mean
as $n\to\infty$.
\end{Lem}

%
\begin{Lem}
\label{THElemmaserial}
For any $\varepsilon>0$, $\sum_{\ell=1}^n \mathrm{ E} [(D^T_{n\ell
})^2 {\mathbb I}[| D^T_{n\ell}| > \varepsilon] ] \to0$
as $n\to\infty$.
\end{Lem}

In the proofs, we use the same notational shortcuts as in the proofs of
Lemmas \ref{firstlemspher}--\ref{THElemmaspher}, that is, we
write $\rho_{st}$, $\Ub_{t}$, $\mathrm{ E}_{\ell}$, $\mathrm{ E}$,
and $\mathrm{
Var}$, instead of $\rho_{n,st}$, $\Ub_{nt}$, $\mathrm{ E}_{n\ell}$,
$\mathrm{
E}_n$ and $\operatorname{ Var}_n$, respectively. For any $r\times s$ matrix
$\Ab
$, we denote as usual by $\operatorname{ vec } \Ab$ the ($rs$)-vector obtained
by stacking the columns of $\Ab$ on top of each other. Recall that we
then have $(\operatorname{ vec } \Ab)'(\operatorname{ vec } \Bb)=\operatorname{
Trace}[\Ab'\Bb]$.

\begin{pf*}{Proof of Lemma~\ref{firstlemserial}}
First note that, for any $s,t\in\{1,\ldots,\ell-1\}$ and any $h,g\in
\{
1,\ldots,H\}$, the $\mathrm{ E}_{\ell-1}$ expectation of
\begin{eqnarray*}
&&\rho_{s-h,\ell-h} \rho_{s,\ell} \rho_{t-g,\ell-g}
\rho_{t,\ell}\\
&&\quad = \bigl( \operatorname{ vec} \bigl(\Ub_s \Ub_{s-h}'
\bigr) \bigr)^{\prime}\operatorname{ vec} \bigl(\Ub_\ell
\Ub_{\ell-h}' \bigr) \bigl( \operatorname{ vec} \bigl(
\Ub_\ell\Ub_{\ell-g}' \bigr) \bigr)^{\prime}
\operatorname{ vec} \bigl(\Ub_t \Ub_{t-g}' \bigr)
\end{eqnarray*}
is given by (see Lemma~\ref{LemmaU}(iv)--(v))
\begin{eqnarray*}
&&\bigl( \operatorname{ vec} \bigl(\Ub_s \Ub_{s-h}'
\bigr) \bigr)^{\prime}\mathrm{ E} \bigl[ \operatorname{ vec} \bigl(\Ub
_\ell \Ub_{\ell-h}' \bigr) \bigl( \operatorname{ vec}
\bigl( \Ub_\ell\Ub_{\ell-g}' \bigr)
\bigr)^{\prime} \bigr] \operatorname{ vec} \bigl(\Ub_t
\Ub_{t-g}' \bigr)
\\
& &\quad = \bigl( \operatorname{ vec} \bigl(\Ub_s \Ub_{s-h}'
\bigr) \bigr)^{\prime}\mathrm{ E} \bigl[ \operatorname{ vec} \bigl(\Ub
_\ell \Ub_{\ell-h}' \bigr) \bigl( \operatorname{ vec}
\bigl( \Ub_\ell\Ub_{\ell-g}' \bigr)
\bigr)^{\prime} \bigr] \operatorname{ vec} \bigl(\Ub_t
\Ub_{t-g}' \bigr)
\\
& &\quad = \frac{1}{p_n^2} \bigl( \operatorname{ vec} \bigl(\Ub_s
\Ub_{s-h}' \bigr) \bigr)^{\prime}\operatorname{ vec} \bigl(
\Ub_t \Ub_{t-g}' \bigr)
\delta_{h,g},
\end{eqnarray*}
where $\delta_{h,g}$ is equal to one if $h=g$ and to zero otherwise.
Therefore, we have that
%
%
\begin{eqnarray}
\label{sigma2serial} \sigma^2_{n\ell} 
& =& \frac{2p_n^2}{H} \sum_{h,g=1}^H
\frac{1}{(n-h)(n-g)} \sum_{s=h+1}^{\ell-1} \sum
_{t=g+1}^{\ell-1} \mathrm{ E}_{\ell-1} [
\rho_{s-h,\ell-h} \rho_{s,\ell} \rho_{t-g,\ell-g}
\rho_{t,\ell} ]
\nonumber
\\
& = & \frac{2}{H} \sum_{h=1}^H
\frac{1}{(n-h)^2} \sum_{s,t=h+1}^{\ell-1} \bigl(
\operatorname{ vec} \bigl(\Ub_s \Ub_{s-h}' \bigr)
\bigr)^{\prime}\operatorname{ vec} \bigl(\Ub_t \Ub_{t-h}'
\bigr)
\\
& = & \frac{2}{H} \sum_{h=1}^H
\frac{1}{(n-h)^2} \sum_{s,t=h+1}^{\ell-1}
\rho_{s-h,t-h} \rho_{s,t} .\nonumber
\end{eqnarray}
From Lemma~\ref{Lemmamoments}(iv), we then obtain
%
%
\begin{equation}
\label{tobeusedserial} \mathrm{ E} \bigl[ \sigma^2_{n\ell} \bigr] =
\frac{2}{H} \sum_{h=1}^H
\frac{(\ell-h-1)_+}{(n-h)^2}\qquad \bigl( = \mathrm{ E} \bigl[ \bigl(D^T_{n\ell}
\bigr)^2 \bigr] \bigr),
\end{equation}
where we let $m_+=\max(m,0)$. This implies that
\[
\mathrm{ E} \Biggl[\sum_{\ell=1}^n
\sigma^2_{n\ell} \Biggr] = \frac{2}{H} \sum
_{h=1}^H \frac{1}{(n-h)^2} \sum
_{\ell=h+2}^n (\ell-h-1) = \frac{1}{H} \sum
_{h=1}^H \frac{n-h-1}{n-h} \to1 ,
\]
as $n\to\infty$.
Using the identity $\operatorname{ Var} [\sum_{h=1}^H Z_h ]\leq H\sum_{h=1}^H \operatorname{ Var} [Z_h ]$, (\ref{sigma2serial}) directly yields
%
%
%
%
\begin{eqnarray}
\label{toplugruns} 
\operatorname{ Var} \Biggl[\sum
_{\ell=1}^n \sigma^2_{n\ell}
\Biggr] & \leq& \frac{4}{H^2} \times H\times\sum
_{h=1}^H \frac{1}{(n-h)^4} \operatorname{ Var} \Biggl[
\sum_{\ell=1}^n \sum
_{s,t=h+1}^{\ell-1} \rho_{s-h,t-h} \rho_{s,t}
\Biggr]
\nonumber
\\[-8pt]
\\[-8pt]
\nonumber
&
\leq& \frac{16}{H(n-H)^4} \sum_{h=1}^H
\operatorname{ Var} \Biggl[ \sum_{\ell=1}^n \sum
_{h+1\leq s<t\leq\ell-1} \rho_{s-h,t-h} \rho_{s,t}
\Biggr] .
\end{eqnarray}
From Lemma~\ref{Lemmamoments}(v), we obtain
\begin{eqnarray*}
&&\operatorname{ Var} \Biggl[ \sum_{\ell=1}^n \sum
_{h+1\leq s<t\leq\ell-1} \rho_{s-h,t-h} \rho_{s,t}
\Biggr]\\
&&\quad \leq\operatorname{ Var} \biggl[ \sum_{h+1\leq s<t\leq n} (n-t)
\rho_{s-h,t-h} \rho_{st} \biggr]
\\
&&\quad  = \frac{1}{p_n^2} \sum_{h+1\leq s<t\leq n}
(n-t)^2 = \frac{1}{p_n^2} \sum_{t=h+2}^n
(t-h-1) (n-t)^2 \\
&&\quad\leq\frac{1}{p_n^2} \sum
_{t=1}^{n-h-1} t (n-h-1)^2
\\
&&\quad  = \frac{1}{p_n^2} \sum_{t=1}^{n-h-1} t
(n-t-h-1)^2 \leq\frac{n^2}{p_n^2} \sum_{t=1}^{n-h-1}
t = \frac{n^2(n-h-1)(n-h)}{2p_n^2}.
\end{eqnarray*}
By plugging this into (\ref{toplugruns}), we obtain
\begin{eqnarray*}
\operatorname{ Var} \Biggl[\sum_{\ell=1}^n
\sigma^2_{n\ell} \Biggr] &\leq& \frac{16}{H(n-H)^4} \sum
_{h=1}^H \frac{n^2(n-h-1)(n-h)}{2p_n^2} \leq
\frac{8n^4}{(n-H)^4p_n^2},
\end{eqnarray*}
which is $\mathrm{o}(1)$ as $n\to\infty$. The result follows.
\end{pf*}

\begin{pf*}{Proof of Lemma~\ref{THElemmaserial}}
Proceeding as in the proof of Lemma~\ref{THElemmaspher}, we obtain
%
%
\begin{eqnarray}
\label{mainserial} \sum_{\ell=1}^n \mathrm{ E}
\bigl[ \bigl(D^T_{n\ell} \bigr)^2 {\mathbb I}
\bigl[\bigl| D^T_{n\ell
}\bigr| > \varepsilon \bigr] \bigr] &\leq&
\frac{1}{\varepsilon} \sum_{\ell=1}^n \sqrt
{\mathrm{ E} \bigl[ \bigl(D^T_{n\ell}
\bigr)^4 \bigr]} \sqrt{\operatorname{ Var} \bigl[D^T_{n\ell}
\bigr]}
\nonumber
\\
&\leq& \frac{\sqrt{2}}{\varepsilon\sqrt{H}} \sum_{h=1}^H
\frac{1}{n-h} \sum_{\ell=h+2}^n \sqrt
{(\ell-h-1) \mathrm{ E} \bigl[ \bigl(D^T_{n\ell}
\bigr)^4 \bigr]}
\\
&\leq& \frac{\sqrt{2}}{\varepsilon(n-H)\sqrt{H}} \sum_{h=1}^H
\sum_{\ell=h+2}^n \sqrt{(\ell-h-1)
\mathrm{ E} \bigl[ \bigl(D^T_{n\ell} \bigr)^4
\bigr]} ,\nonumber
\end{eqnarray}
where we have used the fact that (see (\ref{tobeusedserial}))
\[
\operatorname{ Var} \bigl[D^T_{n\ell} \bigr] \leq\mathrm{ E}
\bigl[ \bigl(D^T_{n\ell} \bigr)^2 \bigr] =
\frac{2}{H} \sum_{h=1}^H
\frac{(\ell-h-1)_+}{(n-h)^2} .
\]
Note that
\begin{eqnarray*}
\label{ok}
&& \mathrm{ E} \Biggl[ \Biggl(\sum_{s=h+1}^{\ell-1}
\rho_{s-h,\ell-h} \rho_{s\ell} \Biggr)^4 \Biggr]
\\
& &\quad = \sum_{s,t,i,j=h+1}^{\ell-1} \mathrm{ E} [
\rho_{s-h,\ell-h} \rho_{t-h,\ell-h} \rho_{i-h,\ell-h}
\rho_{j-h,\ell-h} \rho_{s\ell} \rho_{t\ell}
\rho_{i\ell} \rho_{j \ell} ]
\\
& & \quad= (\ell-h-1)_+ \mathrm{ E} \bigl[ \rho_{1,\ell-h}^4 \bigr]
\mathrm{ E} \bigl[ \rho_{h+1,\ell}^4 \bigr]
\\
& &\qquad{} + 3(\ell-h-1)_+(\ell-h-2)_+ \mathrm{ E} \bigl[ \rho_{1,\ell-h}^2
\bigr] \mathrm{ E} \bigl[ \rho_{2,\ell-h}^2 \bigr] \mathrm{ E}
\bigl[ \rho_{h+1,\ell}^2 \bigr] \mathrm{ E} \bigl[
\rho_{h+2,\ell}^2 \bigr]
\\
& &\quad = \frac{3(\ell-h-1)_+}{p_n^2(p_n+2)^2} + \frac{3(\ell
-h-1)_+(\ell-h-2)_+}{p_n^4} \leq\frac{6(\ell-h-1)_+^2}{p_n^4},
\end{eqnarray*}
which yields
\begin{eqnarray*}
\mathrm{ E} \bigl[ \bigl(D^T_{n\ell} \bigr)^4
\bigr] & \leq& \frac{4p_n^4}{H^2} \times H^{3} \times\sum
_{h=1}^H \frac{1}{(n-h)^4} \mathrm{ E} \Biggl[
\Biggl(\sum_{s=h+1}^{\ell-1} \rho_{s-h,\ell-h}
\rho_{s\ell} \Biggr)^4 \Biggr]
\\
& \leq& \frac{24H}{(n-H)^4} \sum_{h=1}^H (
\ell-h-1)_+^2 .
\end{eqnarray*}
Plugging into (\ref{mainserial}), we conclude that
\begin{eqnarray*}
\sum_{\ell=1}^n \mathrm{ E} \bigl[
\bigl(D^T_{n\ell} \bigr)^2 {\mathbb I} \bigl[\bigl|
D^T_{n\ell
}\bigr| > \varepsilon \bigr] \bigr] & \leq&
\frac{\sqrt{48}}{\varepsilon(n-H)^3} \sum_{h,g=1}^H \sum
_{\ell=h+2}^n \sqrt{(\ell-h-1) (
\ell-g-1)_+^2}
\\
& \leq& \frac{\sqrt{48}H^2}{\varepsilon(n-H)^3} \sum_{\ell=3}^n
(\ell-2)^{3/2}.
\end{eqnarray*}
In view of Lemma~\ref{Lemmaseries}, this is $\mathrm{O}(n^{-1/2})$, which
establishes the result.
\end{pf*}

As an alternative to the tests in (\ref{Tchideux}), one may consider
(see Paindaveine \cite{Pai2009}) the lower-rank multivariate runs
tests rejecting
the null at asymptotic level $\alpha$ whenever
%
%
\[
\tilde{T}_p^{(n)} = \sum_{h=1}^H
(n-h) \bigl( \tilde{\mathbf{r}}(h) \bigr)^2
= \sum
_{h=1}^H \frac{p}{n-h} \sum
_{s,t=h+1}^n \bigl(\Ub_{s-h}'
\Ub_{s} \bigr) \bigl(\Ub_{t-h}'
\Ub_t \bigr) > \chi^2_{H,1-\alpha}, \label{Tchideuxmarden}
\]
where
$
\tilde{\mathbf{r}}(h)
=
\frac{\sqrt{p}}{n-h}
\sum_{t=h+1}^n \Ub_t' \Ub_{t-h}$. Using similar arguments as above,
one may then show that, under the same assumptions as in Theorem~\ref
{maintheoremserial}, the universal $(n,p)$-asymptotic distribution of
$\tilde{T}_{\mathcal{N},p}^{(n)}=(\tilde{T}_n-H)/\sqrt{2H}$
is standard normal.
\end{appendix}

\begin{supplement}
\stitle{Supplement to ``High-dimensional sign tests''}
\slink[doi]{10.3150/15-BEJ710SUPP} 
\sdatatype{.pdf}
\sfilename{BEJ710\_supp.pdf}
\sdescription{The supplement article
contains the proofs of Theorems \ref{maintheoremindep} and \ref{maintheoremsphericity} together with simulation
results related to the sign test for independence. It also provides histograms
from the simulations of Section \ref{simusec1}.}
\end{supplement}

%
%


%
%



\begin{thebibliography}{40}

\bibitem{banerjee2003generative}
%
\begin{binproceedings}[auto:parserefs-M02]
\bauthor{\bsnm{Banerjee},~\bfnm{A.}\binits{A.}},
\bauthor{\bsnm{Dhillon},~\bfnm{I.}\binits{I.}},
\bauthor{\bsnm{Ghosh},~\bfnm{J.}\binits{J.}} \AND
\bauthor{\bsnm{Sra},~\bfnm{S.}\binits{S.}}
(\byear{2003}).
\btitle{Generative model-based clustering of directional data}.
In \bbooktitle{Proceedings of the Ninth ACM SIGKDD International
Conference on Knowledge Discovery and Data Mining}
\bpages{19--28}.
\blocation{New York}:
\bpublisher{ACM}.
\end{binproceedings}
%
\bptok{imsref}%
\endbibitem

\bibitem{Bil1995}
%
\begin{bbook}[mr]
\bauthor{\bsnm{Billingsley},~\bfnm{Patrick}\binits{P.}}
(\byear{1995}).
\btitle{Probability and Measure},
\bedition{3rd} ed.
\bseries{Wiley Series in Probability and Mathematical Statistics}.
\blocation{New York}:
\bpublisher{Wiley}.
\bid{mr={1324786}}
\end{bbook}
%
\bptok{imsref}%
\endbibitem

\bibitem{BroDav1991}
%
\begin{bbook}[mr]
\bauthor{\bsnm{Brockwell},~\bfnm{Peter~J.}\binits{P.J.}} \AND
\bauthor{\bsnm{Davis},~\bfnm{Richard~A.}\binits{R.A.}}
(\byear{1991}).
\btitle{Time Series: Theory and Methods},
\bedition{2nd} ed.
\bseries{Springer Series in Statistics}.
\blocation{New York}:
\bpublisher{Springer}.
\bid{doi={10.1007/978-1-4419-0320-4}, mr={1093459}}
\end{bbook}
%
\bptok{imsref}%
\endbibitem

\bibitem{Caietal2013}
%
\begin{barticle}[mr]
\bauthor{\bsnm{Cai},~\bfnm{Tony}\binits{T.}},
\bauthor{\bsnm{Fan},~\bfnm{Jianqing}\binits{J.}} \AND
\bauthor{\bsnm{Jiang},~\bfnm{Tiefeng}\binits{T.}}
(\byear{2013}).
\btitle{Distributions of angles in random packing on spheres}.
\bjournal{J. Mach. Learn. Res.}
\bvolume{14}
\bpages{1837--1864}.
\bid{issn={1532-4435}, mr={3104497}}
\end{barticle}
%
\bptok{imsref}%
\endbibitem

\bibitem{CaJia12}
%
\begin{barticle}[mr]
\bauthor{\bsnm{Cai},~\bfnm{T.~Tony}\binits{T.T.}} \AND
\bauthor{\bsnm{Jiang},~\bfnm{Tiefeng}\binits{T.}}
(\byear{2012}).
\btitle{Phase transition in limiting distributions of coherence of
high-dimensional random matrices}.
\bjournal{J. Multivariate Anal.}
\bvolume{107}
\bpages{24--39}.
\bid{doi={10.1016/j.jmva.2011.11.008}, issn={0047-259X}, mr={2890430}}
\end{barticle}
%
\bptok{imsref}%
\endbibitem

\bibitem{Cheetal2010}
%
\begin{barticle}[mr]
\bauthor{\bsnm{Chen},~\bfnm{Song~Xi}\binits{S.X.}},
\bauthor{\bsnm{Zhang},~\bfnm{Li-Xin}\binits{L.-X.}} \AND
\bauthor{\bsnm{Zhong},~\bfnm{Ping-Shou}\binits{P.-S.}}
(\byear{2010}).
\btitle{Tests for high-dimensional covariance matrices}.
\bjournal{J.~Amer. Statist. Assoc.}
\bvolume{105}
\bpages{810--819}.
\bid{doi={10.1198/jasa.2010.tm09560}, issn={0162-1459}, mr={2724863}}
\end{barticle}
%
\bptok{imsref}%
\endbibitem

\bibitem{Chi1991}
%
\begin{barticle}[mr]
\bauthor{\bsnm{Chikuse},~\bfnm{Yasuko}\binits{Y.}}
(\byear{1991}).
\btitle{High-dimensional limit theorems and matrix decompositions on
the {S}tiefel manifold}.
\bjournal{J. Multivariate Anal.}
\bvolume{36}
\bpages{145--162}.
\bid{doi={10.1016/0047-259X(91)90054-6}, issn={0047-259X}, mr={1096663}}
\end{barticle}
%
\bptok{imsref}%
\endbibitem

\bibitem{Dry2005}
%
\begin{barticle}[mr]
\bauthor{\bsnm{Dryden},~\bfnm{Ian~L.}\binits{I.L.}}
(\byear{2005}).
\btitle{Statistical analysis on high-dimensional spheres and shape spaces}.
\bjournal{Ann. Statist.}
\bvolume{33}
\bpages{1643--1665}.
\bid{doi={10.1214/009053605000000264}, issn={0090-5364}, mr={2166558}}
\end{barticle}
%
\bptok{imsref}%
\endbibitem

\bibitem{Dufetal1998}
%
\begin{barticle}[mr]
\bauthor{\bsnm{Dufour},~\bfnm{Jean-Marie}\binits{J.-M.}},
\bauthor{\bsnm{Hallin},~\bfnm{Marc}\binits{M.}} \AND
\bauthor{\bsnm{Mizera},~\bfnm{Ivan}\binits{I.}}
(\byear{1998}).
\btitle{Generalized runs tests for heteroscedastic time series}.
\bjournal{J. Nonparametr. Stat.}
\bvolume{9}
\bpages{39--86}.
\bid{doi={10.1080/10485259808832735}, issn={1048-5252}, mr={1623129}}
\end{barticle}
%
\bptok{imsref}%
\endbibitem



\bibitem{Dum98}
%
\begin{barticle}[mr]
\bauthor{\bsnm{D{\"u}mbgen},~\bfnm{Lutz}\binits{L.}}
(\byear{1998}).
\btitle{On {T}yler's {$M$}-functional of scatter in high dimension}.
\bjournal{Ann. Inst. Statist. Math.}
\bvolume{50}
\bpages{471--491}.
\bid{doi={10.1023/A:1003573311481}, issn={0020-3157}, mr={1664575}}
\end{barticle}
%
\bptok{imsref}%
\endbibitem

\bibitem{fra2004}
%
\begin{bmisc}[auto]
\bauthor{\bsnm{Frahm},~\bfnm{Gabriel}\binits{G.}}
(\byear{2004}).
\bhowpublished{Generalized elliptical distributions: Theory
and applications. Ph.D. thesis, Universit\"{a}t zu K\"{o}ln.}
\end{bmisc}
\bptok{imsref}%
\endbibitem


\bibitem{frajae2007}
%
\begin{bmisc}[auto]
\bauthor{\bsnm{Frahm},~\bfnm{Gabriel}\binits{G.}} \AND
\bauthor{\bsnm{Jaekel},~\bfnm{Uwe}\binits{U.}}
(\byear{2010}).
\bhowpublished{Tyler's $M$-estimator, random
matrix theory, and generalized elliptical distributions with applications
to finance. Discussion paper, Department of Economic and Social Statistics,
University of Cologne, Germany.}
\end{bmisc}
\bptok{imsref}%
\endbibitem

\bibitem{frajae2010}
%
\begin{barticle}[mr]
\bauthor{\bsnm{Frahm},~\bfnm{Gabriel}\binits{G.}} \AND
\bauthor{\bsnm{Jaekel},~\bfnm{Uwe}\binits{U.}}
(\byear{2010}).
\btitle{A generalization of {T}yler's {$M$}-estimators to the case of
incomplete data}.
\bjournal{Comput. Statist. Data Anal.}
\bvolume{54}
\bpages{374--393}.
\bid{doi={10.1016/j.csda.2009.08.019}, issn={0167-9473}, mr={2756433}}
\end{barticle}
%
\bptok{imsref}%
\endbibitem

\bibitem{Halletal2005}
%
\begin{barticle}[mr]
\bauthor{\bsnm{Hall},~\bfnm{Peter}\binits{P.}},
\bauthor{\bsnm{Marron},~\bfnm{J.~S.}\binits{J.S.}} \AND
\bauthor{\bsnm{Neeman},~\bfnm{Amnon}\binits{A.}}
(\byear{2005}).
\btitle{Geometric representation of high dimension, low sample size data}.
\bjournal{J. R. Stat. Soc. Ser. B. Stat. Methodol.}
\bvolume{67}
\bpages{427--444}.
\bid{doi={10.1111/j.1467-9868.2005.00510.x}, issn={1369-7412}, mr={2155347}}
\end{barticle}
%
\bptok{imsref}%
\endbibitem

\bibitem{HP02}
%
\begin{barticle}[mr]
\bauthor{\bsnm{Hallin},~\bfnm{Marc}\binits{M.}} \AND
\bauthor{\bsnm{Paindaveine},~\bfnm{Davy}\binits{D.}}
(\byear{2002}).
\btitle{Optimal tests for multivariate location based on
interdirections and pseudo-{M}ahalanobis ranks}.
\bjournal{Ann. Statist.}
\bvolume{30}
\bpages{1103--1133}.
\bid{doi={10.1214/aos/1031689019}, issn={0090-5364}, mr={1926170}}
\end{barticle}
%
\bptok{imsref}%
\endbibitem

\bibitem{HP06}
%
\begin{barticle}[mr]
\bauthor{\bsnm{Hallin},~\bfnm{Marc}\binits{M.}} \AND
\bauthor{\bsnm{Paindaveine},~\bfnm{Davy}\binits{D.}}
(\byear{2006}).
\btitle{Semiparametrically efficient rank-based inference for shape.
I. {O}ptimal rank-based tests for sphericity}.
\bjournal{Ann. Statist.}
\bvolume{34}
\bpages{2707--2756}.
\bid{doi={10.1214/009053606000000731}, issn={0090-5364}, mr={2329465}}
\end{barticle}
%
\bptok{imsref}%
\endbibitem

\bibitem{HPV10}
%
\begin{barticle}[mr]
\bauthor{\bsnm{Hallin},~\bfnm{Marc}\binits{M.}},
\bauthor{\bsnm{Paindaveine},~\bfnm{Davy}\binits{D.}} \AND
\bauthor{\bsnm{Verdebout},~\bfnm{Thomas}\binits{T.}}
(\byear{2010}).
\btitle{Optimal rank-based testing for principal components}.
\bjournal{Ann. Statist.}
\bvolume{38}
\bpages{3245--3299}.
\bid{doi={10.1214/10-AOS810}, issn={0090-5364}, mr={2766852}}
\end{barticle}
%
\bptok{imsref}%
\endbibitem

\bibitem{JiaYan2013}
%
\begin{barticle}[mr]
\bauthor{\bsnm{Jiang},~\bfnm{Tiefeng}\binits{T.}} \AND
\bauthor{\bsnm{Yang},~\bfnm{Fan}\binits{F.}}
(\byear{2013}).
\btitle{Central limit theorems for classical likelihood ratio tests
for high-dimensional normal distributions}.
\bjournal{Ann. Statist.}
\bvolume{41}
\bpages{2029--2074}.
\bid{doi={10.1214/13-AOS1134}, issn={0090-5364}, mr={3127857}}
\bptnote{check pages}%
\end{barticle}
%
\bptok{imsref}%
\endbibitem

\bibitem{Johetal1995}
%
\begin{bbook}[mr]
\bauthor{\bsnm{Johnson},~\bfnm{Norman~L.}\binits{N.L.}},
\bauthor{\bsnm{Kotz},~\bfnm{Samuel}\binits{S.}} \AND
\bauthor{\bsnm{Balakrishnan},~\bfnm{N.}\binits{N.}}
(\byear{1995}).
\btitle{Continuous Univariate Distributions. {V}ol. 2},
\bedition{2nd} ed.
\bseries{Wiley Series in Probability and Mathematical Statistics:
Applied Probability and Statistics}.
\blocation{New York}:
\bpublisher{Wiley}.
\bid{mr={1326603}}
\end{bbook}
%
\bptok{imsref}%
\endbibitem

\bibitem{LedWol2002}
%
\begin{barticle}[mr]
\bauthor{\bsnm{Ledoit},~\bfnm{Olivier}\binits{O.}} \AND
\bauthor{\bsnm{Wolf},~\bfnm{Michael}\binits{M.}}
(\byear{2002}).
\btitle{Some hypothesis tests for the covariance matrix when the
dimension is large compared to the sample size}.
\bjournal{Ann. Statist.}
\bvolume{30}
\bpages{1081--1102}.
\bid{doi={10.1214/aos/1031689018}, issn={0090-5364}, mr={1926169}}
\end{barticle}
%
\bptok{imsref}%
\endbibitem

\bibitem{Leyetal2013}
%
\begin{barticle}[mr]
\bauthor{\bsnm{Ley},~\bfnm{Christophe}\binits{C.}},
\bauthor{\bsnm{Swan},~\bfnm{Yvik}\binits{Y.}},
\bauthor{\bsnm{Thiam},~\bfnm{Baba}\binits{B.}} \AND
\bauthor{\bsnm{Verdebout},~\bfnm{Thomas}\binits{T.}}
(\byear{2013}).
\btitle{Optimal {$R$}-estimation of a spherical location}.
\bjournal{Statist. Sinica}
\bvolume{23}
\bpages{305--332}.
\bid{issn={1017-0405}, mr={3076169}}
\bptnote{check pages}%
\end{barticle}
%
\bptok{imsref}%
\endbibitem

\bibitem{li1981}
%
\begin{barticle}[mr]
\bauthor{\bsnm{Li},~\bfnm{W.~K.}\binits{W.K.}} \AND
\bauthor{\bsnm{McLeod},~\bfnm{A.~I.}\binits{A.I.}}
(\byear{1981}).
\btitle{Distribution of the residual autocorrelations in multivariate
ARMA time series models}.
\bjournal{J. R. Stat. Soc. Ser. B. Stat. Methodol.}
\bvolume{43}
\bpages{231--239}.
\bid{issn={0035-9246}, mr={0626770}}
\end{barticle}
%
\bptok{imsref}%
\endbibitem

\bibitem{ljung78}
%
\begin{barticle}[auto:parserefs-M02]
\bauthor{\bsnm{Ljung},~\bfnm{G.~M.}\binits{G.M.}} \AND
\bauthor{\bsnm{Box},~\bfnm{G.~E.}\binits{G.E.}}
(\byear{1978}).
\btitle{On a measure of lack of fit in time series models}.
\bjournal{Biometrika}
\bvolume{65}
\bpages{297--303}.
\end{barticle}
%
\bptok{imsref}%
\endbibitem

\bibitem{Lut2005}
%
\begin{bbook}[mr]
\bauthor{\bsnm{L{\"u}tkepohl},~\bfnm{Helmut}\binits{H.}}
(\byear{2005}).
\btitle{New Introduction to Multiple Time Series Analysis}.
\blocation{Berlin}:
\bpublisher{Springer}.
\bid{doi={10.1007/978-3-540-27752-1}, mr={2172368}}
\end{bbook}
%
\bptok{imsref}%
\endbibitem

\bibitem{MarJup2000}
%
\begin{bbook}[mr]
\bauthor{\bsnm{Mardia},~\bfnm{Kanti~V.}\binits{K.V.}} \AND
\bauthor{\bsnm{Jupp},~\bfnm{Peter~E.}\binits{P.E.}}
(\byear{2000}).
\btitle{Directional Statistics}.
\bseries{Wiley Series in Probability and Statistics}.
\blocation{Chichester}:
\bpublisher{Wiley}.
\bid{mr={1828667}}
\end{bbook}
%
\bptok{imsref}%
\endbibitem

\bibitem{mooj95}
%
\begin{barticle}[mr]
\bauthor{\bsnm{M{\"o}tt{\"o}nen},~\bfnm{Jyrki}\binits{J.}} \AND
\bauthor{\bsnm{Oja},~\bfnm{Hannu}\binits{H.}}
(\byear{1995}).
\btitle{Multivariate spatial sign and rank methods}.
\bjournal{J. Nonparametr. Stat.}
\bvolume{5}
\bpages{201--213}.
\bid{doi={10.1080/10485259508832643}, issn={1048-5252}, mr={1346895}}
\end{barticle}
%
\bptok{imsref}%
\endbibitem

\bibitem{Mui1982}
%
\begin{bbook}[mr]
\bauthor{\bsnm{Muirhead},~\bfnm{Robb~J.}\binits{R.J.}}
(\byear{1982}).
\btitle{Aspects of Multivariate Statistical Theory}.
\blocation{New York}:
\bpublisher{Wiley}.
\bid{mr={0652932}}
\end{bbook}
%
\bptok{imsref}%
\endbibitem

\bibitem{Oja2010}
%
\begin{bbook}[mr]
\bauthor{\bsnm{Oja},~\bfnm{Hannu}\binits{H.}}
(\byear{2010}).
\btitle{Multivariate Nonparametric Methods with R: An Approach Based
on Spatial Signs and Ranks}.
\bseries{Lecture Notes in Statistics}
\bvolume{199}.
\blocation{New York}:
\bpublisher{Springer}.
\bid{doi={10.1007/978-1-4419-0468-3}, mr={2598854}}
\end{bbook}
%
\bptok{imsref}%
\endbibitem

\bibitem{Onaetal2013}
%
\begin{barticle}[mr]
\bauthor{\bsnm{Onatski},~\bfnm{Alexei}\binits{A.}},
\bauthor{\bsnm{Moreira},~\bfnm{Marcelo~J.}\binits{M.J.}} \AND
\bauthor{\bsnm{Hallin},~\bfnm{Marc}\binits{M.}}
(\byear{2013}).
\btitle{Asymptotic power of sphericity tests for high-dimensional data}.
\bjournal{Ann. Statist.}
\bvolume{41}
\bpages{1204--1231}.
\bid{doi={10.1214/13-AOS1100}, issn={0090-5364}, mr={3113808}}
\end{barticle}
%
\bptok{imsref}%
\endbibitem

\bibitem{Pai2009}
%
\begin{barticle}[mr]
\bauthor{\bsnm{Paindaveine},~\bfnm{Davy}\binits{D.}}
(\byear{2009}).
\btitle{On multivariate runs tests for randomness}.
\bjournal{J. Amer. Statist. Assoc.}
\bvolume{104}
\bpages{1525--1538}.
\bid{doi={10.1198/jasa.2009.tm09047}, issn={0162-1459}, mr={2597001}}
\end{barticle}
%
\bptok{imsref}%
\endbibitem

\bibitem{PaiVer2013}
%
\begin{bincollection}[auto:parserefs-M02]
\bauthor{\bsnm{Paindaveine},~\bfnm{D.}\binits{D.}} \AND
\bauthor{\bsnm{Verdebout},~\bfnm{T.}\binits{T.}}
(\byear{2014}).
\btitle{Optimal rank-based tests for the location parameter of a
rotationally symmetric distribution on the hypersphere}.
In \bbooktitle{Mathematical Statistics and Limit Theorems: Festschrift
in Honor of Paul Deheuvels}
(\beditor{\bfnm{M.}\binits{M.}~\bsnm{Hallin}},
\beditor{\bfnm{D.}\binits{D.}~\bsnm{Mason}},
\beditor{\bfnm{D.}\binits{D.}~\bsnm{Pfeifer}} \AND
\beditor{\bfnm{J.}\binits{J.}~\bsnm{Steinebach}}, eds.)
\bpages{249--270}.
\blocation{Berlin}:
\bpublisher{Springer}.
\end{bincollection}
%
\bptok{imsref}%
\endbibitem

\bibitem{PaiVer2013c}
%
\begin{bmisc}[author]
\bauthor{\bsnm{Paindaveine},~\bfnm{D.}\binits{D.}} \AND
\bauthor{\bsnm{Verdebout},~\bfnm{T.}\binits{T.}}
(\byear{2015}).
\bhowpublished{Supplement to ``On high-dimensional sign tests.''
DOI:\doiurl{10.3150/15-BEJ710SUPP}}.
\bptok{imsref}%
\end{bmisc}
%
\endbibitem%

\bibitem{Ran89}
%
\begin{barticle}[mr]
\bauthor{\bsnm{Randles},~\bfnm{Ronald~H.}\binits{R.H.}}
(\byear{1989}).
\btitle{A distribution-free multivariate sign test based on interdirections}.
\bjournal{J. Amer. Statist. Assoc.}
\bvolume{84}
\bpages{1045--1050}.
\bid{issn={0162-1459}, mr={1134492}}
\end{barticle}
%
\bptok{imsref}%
\endbibitem

\bibitem{Ran00}
%
\begin{barticle}[mr]
\bauthor{\bsnm{Randles},~\bfnm{Ronald~H.}\binits{R.H.}}
(\byear{2000}).
\btitle{A simpler, affine-invariant, multivariate, distribution-free
sign test}.
\bjournal{J. Amer. Statist. Assoc.}
\bvolume{95}
\bpages{1263--1268}.
\bid{doi={10.2307/2669766}, issn={0162-1459}, mr={1792189}}
\end{barticle}
%
\bptok{imsref}%
\endbibitem

\bibitem{Ray1919}
%
\begin{barticle}[auto:parserefs-M02]
\bauthor{\bsnm{Rayleigh},~\bfnm{L.}\binits{L.}}
(\byear{1919}).
\btitle{On the problem of random vibrations and random flights in one,
two and three dimensions}.
\bjournal{Phil. Mag.}
\bvolume{37}
\bpages{321--346}.
\end{barticle}
%
\bptok{imsref}%
\endbibitem

\bibitem{Siretal2009}
%
\begin{barticle}[mr]
\bauthor{\bsnm{Sirki{\"a}},~\bfnm{Seija}\binits{S.}},
\bauthor{\bsnm{Taskinen},~\bfnm{Sara}\binits{S.}},
\bauthor{\bsnm{Oja},~\bfnm{Hannu}\binits{H.}} \AND
\bauthor{\bsnm{Tyler},~\bfnm{David~E.}\binits{D.E.}}
(\byear{2009}).
\btitle{Tests and estimates of shape based on spatial signs and ranks}.
\bjournal{J. Nonparametr. Stat.}
\bvolume{21}
\bpages{155--176}.
\bid{doi={10.1080/10485250802495691}, issn={1048-5252}, mr={2488152}}
\end{barticle}
%
\bptok{imsref}%
\endbibitem

\bibitem{SriFuj2006}
%
\begin{barticle}[mr]
\bauthor{\bsnm{Srivastava},~\bfnm{Muni~S.}\binits{M.S.}} \AND
\bauthor{\bsnm{Fujikoshi},~\bfnm{Yasunori}\binits{Y.}}
(\byear{2006}).
\btitle{Multivariate analysis of variance with fewer observations than
the dimension}.
\bjournal{J. Multivariate Anal.}
\bvolume{97}
\bpages{1927--1940}.
\bid{doi={10.1016/j.jmva.2005.08.010}, issn={0047-259X}, mr={2301265}}
\end{barticle}
%
\bptok{imsref}%
\endbibitem

\bibitem{SriKub2013}
%
\begin{barticle}[mr]
\bauthor{\bsnm{Srivastava},~\bfnm{Muni~S.}\binits{M.S.}} \AND
\bauthor{\bsnm{Kubokawa},~\bfnm{Tatsuya}\binits{T.}}
(\byear{2013}).
\btitle{Tests for multivariate analysis of variance in high dimension
under non-normality}.
\bjournal{J. Multivariate Anal.}
\bvolume{115}
\bpages{204--216}.
\bid{doi={10.1016/j.jmva.2012.10.011}, issn={0047-259X}, mr={3004555}}
\end{barticle}
%
\bptok{imsref}%
\endbibitem

\bibitem{Tasetal2003}
%
\begin{barticle}[mr]
\bauthor{\bsnm{Taskinen},~\bfnm{Sara}\binits{S.}},
\bauthor{\bsnm{Kankainen},~\bfnm{Annaliisa}\binits{A.}} \AND
\bauthor{\bsnm{Oja},~\bfnm{Hannu}\binits{H.}}
(\byear{2003}).
\btitle{Sign test of independence between two random vectors}.
\bjournal{Statist. Probab. Lett.}
\bvolume{62}
\bpages{9--21}.
\bid{doi={10.1016/S0167-7152(02)00399-1}, issn={0167-7152}, mr={1965367}}
\end{barticle}
%
\bptok{imsref}%
\endbibitem

\bibitem{Tasetal2005}
%
\begin{barticle}[mr]
\bauthor{\bsnm{Taskinen},~\bfnm{Sara}\binits{S.}},
\bauthor{\bsnm{Oja},~\bfnm{Hannu}\binits{H.}} \AND
\bauthor{\bsnm{Randles},~\bfnm{Ronald~H.}\binits{R.H.}}
(\byear{2005}).
\btitle{Multivariate nonparametric tests of independence}.
\bjournal{J.~Amer. Statist. Assoc.}
\bvolume{100}
\bpages{916--925}.
\bid{doi={10.1198/016214505000000097}, issn={0162-1459}, mr={2201019}}
\end{barticle}
%
\bptok{imsref}%
\endbibitem

\bibitem{Tyl1987}
%
\begin{barticle}[mr]
\bauthor{\bsnm{Tyler},~\bfnm{David~E.}\binits{D.E.}}
(\byear{1987}).
\btitle{A distribution-free {$M$}-estimator of multivariate scatter}.
\bjournal{Ann. Statist.}
\bvolume{15}
\bpages{234--251}.
\bid{doi={10.1214/aos/1176350263}, issn={0090-5364}, mr={0885734}}
\end{barticle}
%
\bptok{imsref}%
\endbibitem

\bibitem{Zouetal2013}
%
\begin{barticle}[mr]
\bauthor{\bsnm{Zou},~\bfnm{Changliang}\binits{C.}},
\bauthor{\bsnm{Peng},~\bfnm{Liuhua}\binits{L.}},
\bauthor{\bsnm{Feng},~\bfnm{Long}\binits{L.}} \AND
\bauthor{\bsnm{Wang},~\bfnm{Zhaojun}\binits{Z.}}
(\byear{2014}).
\btitle{Multivariate sign-based high-dimensional tests for sphericity}.
\bjournal{Biometrika}
\bvolume{101}
\bpages{229--236}.
\bid{doi={10.1093/biomet/ast040}, issn={0006-3444}, mr={3180668}}
\end{barticle}
%
\bptok{imsref}%
\endbibitem
\end{thebibliography}
%
%




\printhistory
\end{document}